\documentclass[12pt]{article}
\usepackage{amsmath}
\usepackage{amssymb}
\usepackage{amsfonts}
\usepackage{eucal}
\usepackage[usenames]{color}
\usepackage{graphicx}
\usepackage[symbol]{footmisc}
\numberwithin{equation}{section}
\newtheorem{Theorem}{Theorem}[section]
\newtheorem{Proposition}[Theorem]{Proposition}

\newtheorem{Remark}[Theorem]{Remark}
\usepackage{upgreek}

\newcommand{\ep}[1]{{\textcolor{black}{#1}}}

\begin{document}
\title{On the Fractional Riemann-Liouville Integral of Gauss-Markov
processes and applications}

\author{Mario Abundo\thanks{Dipartimento di Matematica, Universit\`a
``Tor Vergata'', via della Ricerca Scientifica, I-00133 Roma,
Italy. E-mail: \texttt{abundo@mat.uniroma2.it}} 
\and Enrica Pirozzi
\footnote{Dipartimento di Matematica e Applicazioni, Universit\`a ``Federico II'', via
Cintia, complesso Monte S. Angelo, I-80126 Napoli, Italy.
E-mail: \texttt{enrica.pirozzi@unina.it}}
}

\date{}
\maketitle

\begin{abstract}
\noindent \ep{
We investigate the stochastic processes obtained as the fractional Riemann-Liouville integral of order $\alpha \in (0,1)$ of  Gauss-Markov processes. The general expressions of the mean, variance and covariance functions are given. Due to the central rule, for the fractional integral of standard Brownian motion and of the non-stationary/stationary Ornstein-Uhlenbeck processes, the covariance functions are carried out in closed-form. In order to clarify how the fractional order parameter $\alpha$ affects these functions, their numerical evaluations are shown and compared also with those of the corresponding processes obtained by ordinary Riemann integral.  The results are useful for fractional neuronal models with long range memory dynamics and involving correlated input processes. The simulation of these fractional integrated processes can be performed starting from the obtained covariance functions. A suitable neuronal model is proposed. Graphical comparisons are provided and discussed.}
\end{abstract}

\noindent {\bf Keywords:} Fractional integrals; Covariance function; Neuronal models; Simulation \\
{\bf Mathematics Subject Classification:} 60G15, 26A33, 65C20.

\section{Introduction}

\ep{Although fractional calculus has the same ancient origins of the classical calculus, it has become of extreme interest in last decades for the researches, in many and different science fields. More recently, the theory and applications of fractional derivatives and fractional integrals (\cite{Sam93}) have been  
extensively developed  by pure  and applied mathematicians. Among a lot of motivations, we can claim that the researcher's community realized that fractional differential equations and related fractional integral solutions provide a natural framework for the description and the study of real phenomena such as, for instance, those occurring in biology, in ecology and in neuroscience (see, for instance, \cite{Bal09},\cite{Hil00},\cite{Lak09} and references therein).} 

\ep{ Firstly, the integrals of stochastic processes constitute a mathematical subject fascinating especially probability theorists (see, for instance, \cite{Cui17}, \cite{Ste95} and references therein); successively a wide audience was also attracted. Indeed, the need to describe complex phenomena, whose time evolution is affected by the history of own behavior (\cite{Fer12}) and by several memory effects of different nature (\cite{Jaf17}), requires to design more elaborate mathematical models relying on integrals over time of stochastic processes, and some their modifications. Several aspects, particular cases and specific models can be found in \cite{abundo:smj15}-\cite{abupir:physicaA17}. In particular, in \cite{abupir:physicaA17} we considered  Gauss-Markov (GM) processes and their time integrals. Due to their easy mathematical handling, such kind of processes are suitable for modeling purposes, in a wide range of applications as, for instance,  in computational neuroscience
(see e.g. \cite{tou:aap08} and references therein), in finance
mathematics (see e.g. \cite{bar:03}), in queueing theory and other applied sciences (see e.g. the discussion in
\cite{abundo:smj15}, \cite{abundo:smj13}). 
Moreover, as a specific example of application of integrated GM processes, we can refer to the
context of neuronal modeling; indeed, some dynamics have been studied
by introducing the so-called colored noise (i.e., a correlated GM process) in  neuronal
stochastic models, in place of the classical white noise (see e.g.
\cite{Baz03}, \cite{Kim14}, \cite{Kob09}, \cite{ep:biocyb17}).
This kind of models rely on  stochastic processes which are the
integrals over time of  an Ornstein-Uhlenbeck (OU) process, or more
generally, of a GM process.
Integrated GM processes and their first-passage times (FPT) were
studied first in \cite{abundo:smj15}, \cite{abundo:smj13}, and
then in \cite{abupir:physicaA17}.}

\ep{With the aim to specialize this topic,  in this paper we consider its extension to the fractional integration. Specifically, we investigate the  stochastic processes obtained as Riemann-Liouville (RL) integral of GM processes. The motivation of this research is not only to shed light on such processes, investigating their properties, useful to consolidate their  mathematical setting, but also to explore their powerful skills to model all phenomena that preserve their memory and are resultant of other dynamics evolving on different time-scales.  The fractional integral reveals really suitable to specialize such kind of models; indeed, this mathematical tool plays the key rule for tuning the scale of the time, in such a way new features of the dynamical systems emerge and can be studied. }

\ep {In the framework of
neuronal dynamics, in \cite{ep:biocyb17},
fractional stochastic models are introduced for
preserving the memory of the neuronal membrane evolution; in such
description, the parameter that affects the firing activity is the
fractional order of the involved derivative. Simulation procedures are a first step in the investigation of such kind of models (see, for instance, \cite{Asc18},\cite{ep:biocyb17},\cite{Tek14}).
The aim of  giving a further contribution, that can be useful in the theoretical and simulation approaches to neuronal modeling,  led us
to consider this type of variable range memory process, obtained as the fractional integral over time
of  GM processes. An example of the mathematical extension of a classical model to the fractional one has been recently done, for instance, in a different context such as that of queueing systems in \cite{Asc18bis}}. \par
\ep{In this article, motivated by  all above
considerations, we  study as replacing the ordinary Riemann integral with the fractional RL
integral of order $\alpha \in (0,1)$ affects the behavior, when $\alpha$ varies, of an integrated GM process. 
Due to their importance among GM processes, we mainly focus the attention on the fractional integral of the Brownian motion and  OU processes. In particular, we consider the cases of stationary and non-stationary OU processes, because these two processes allow to design different neuronal models in which stimuli of different nature (exogenous and endogenous stimuli, respectively) can be included (\cite{Asc19}).
A different approach has been considered in \cite{Sit97}, where it was studied 
the power spectral density of RL fractional BM, in order to establish a 
fractional power law of the form $1/ f^ \alpha $; note that, for instance,
this feature is observed in ordinary BM and its time integral.
Moreover, in the present article, we do not even consider the 
fractional OU process as in \cite{Lim15}, in which it is defined as the Riemann-Liuoville of a tempered OU process (for some details on tempered fractional calculus see, e.g., \cite{Mee14} and references therein). Further investigations on such processes and comparisons with the processes here analyzed will be part of  our future work. 
%
}
\par
\ep{Furthermore, as an example for possible applications, we consider the model for the neuronal activity  based on the following coupled differential equations, for $t\geq 0$:}
\ep{
\begin{eqnarray}\label{1}
{\cal D}^\alpha V(t)&=&\frac{g_L}{C_m}V_L+\frac{{{\upeta}(t)}}{C_m}, \hspace{2.2cm} V(0)=V_0\\
d{\upeta(t)}&=& - \frac{{\upeta(t)}-{I(t)}}{\tau}dt+\frac{\varsigma}{\tau}dB(t), \ \upeta(0)=\upeta_0. \label{2}
\end{eqnarray}}
\par\noindent
\ep{{{In the \eqref{1} the derivative ${\cal D}^\alpha$ stands for the Caputo fractional derivative (specified in the next section); moreover, $\upeta(t)$ is in place of the white noise $dB(t)$ as usual in the stochastic differential equation (SDE) of a Leaky Integrate-and-Fire (LIF) neuronal model (see, for example, \cite{Bur06}).}}
The colored noise process $\upeta(t)$ is the correlated process obeying to the SDE \eqref{2} and it is the input for the equation \eqref{1}.
The stochastic process $V(t)$ represents the voltage of the neuronal membrane, whereas the other parameters and functions are: $C_m$ the membrane capacitance, $g_L$  the leak conductance, $V_L$  the resting (equilibrium) level potential, $I(t)$ the synaptic current (deterministic function),  ${\tau}$ is the correlation time of $\upeta(t)$ and $W(t)$ the noise (a standard Brownian motion). 
The initial values $V_0$ and $\upeta_0$ can be specified constants or random variables, defining the above model as stationary or non-stationary one, respectively.
 $\upeta(t)$ is a time-non-homogeneous GM process (OU-type). Finally, the solution process $V(t)$ of the equation \eqref{1} belongs to the class of the fractional RL integrals of the $\upeta(t)$ GM process, that will be studied in this paper.}

\ep{Note that, beyond the colored noise, the novelty of the above neuronal model, respect to the classical ones, is in the use of the fractional derivative in place of the integer one. The biophysical motivation is  to describe a neuronal activity (by equation \eqref{1}) as a perfect integrator (without leakage) of the whole evolution of the input process $\upeta(t)$ from an initial time until to the current time, but on a (more or less) finer time scale that can be regulated by choosing the fractional order of integration suitably adherent to the neuro-physiological evidences. Indeed,  such a model can be useful, for instance, in the investigation and simulation of synchronous/asynchronous communications in networks of neurons (\cite{Tam18}). }


\ep{We remark again that here we study the processes achieved applying the fractional RL integral to GM processes, and we investigate how the variation of the fractional order parameter affects the main features of these processes, such as  mean, variance, covariance and paths. Furthermore, the mathematical results are shown, compared and discussed also in graphical way in  figures carried out by  numerical and simulation procedures. Finally, some indications about their usefulness for neuronal modelling are provided.}

\ep{ The paper is organized as follows. Section 2 contains mathematical formulation and details of the fractional RL integral of the BM (FIBM); plots of variance and covariance functions are provided by means of numerical evaluations; sample paths are also simulated for different values of $\alpha$ in order to show the qualitative behavior of the process.  Our main result is in Section 3: the fractional RL integral of a GM process is defined and the covariance function is evaluated. Then, in Section 4, for the RL fractional integral of OU and stationary OU processes, the covariance functions are calculated and numerically evaluated carrying out comparisons between the fractional integral and the ordinary integral of these specific GM processes. Finally, in Section 5 we report some graphical comparisons and  concluding remarks.}

\section{Mathematical formulation and preliminary results}
\ep{At first, we start recalling some well-known definitions and introducing the object of our main interest.}
  
Let be $\alpha \in (0,1);$ if $f(t)$ is a \ep{real-valued differentiable function on $\mathbb R$}, we recall  that the
{\it Caputo fractional derivative} of $f$ of order $\alpha$ is defined by (see \cite{Cap67}):
\begin{equation} \label{fractderivative}
{\cal D} ^\alpha  f(t) = \frac {1 } {\Gamma (1- \alpha) } \int _0 ^t \frac {f'(s) } {(t-s)^\alpha } \ ds ,
\end{equation}
where $f'$ denotes the  ordinary derivative of $f.$
\par\noindent
If $f$ is a continuous function,
its {\it fractional RL integral} of order $\alpha$ is defined  by (see \cite{Deb03}):
\begin{equation} \label{fractintegral}
{\cal I }^ \alpha  (f)(t) = \frac {1 } {\Gamma (\alpha) } \int _0 ^t (t-s)^{ \alpha -1} f(s) ds ,
\end{equation}
where $\Gamma $ is the Gamma Euler function, i.e. $\Gamma (z)=
\int _0 ^ { + \infty } t^{z-1} e ^{-t} dt \ , \ z >0.$ Notice
that, taking the limit for $\alpha \rightarrow 1 ^-,$
\eqref{fractderivative} provides the ordinary derivative of $f,$
while \eqref{fractintegral} gives the ordinary Riemann integral of
$f;$ moreover, by convention, ${\cal D} ^0  f(t)= f(t)-f(0)$ and ${\cal I }^ 0  (f)(t) = f(t)$
(for properties of Caputo fractional derivative and
fractional RL integral, see e.g. \ep{ \cite{Ish05},\cite{Li11},\cite{Pod99},\cite{Sam93}}).
\par Now, we recall the definition of GM process:
let $m(t), \ h_1(t), \ h_2(t)$ be continuous functions of $t \ge
0$  which are $C^1$ in $(0, + \infty )$ and such that $h_2(t) \neq
0 \; \forall t\geq0,$ and let $ r(t)= h_1(t) / h_2(t) $ be a non-negative,
differentiable function, with $r'(t)>0$ for $t >0,$ and $r (0)=r_0
\ge 0 .$ \par\noindent 
If, \ep{for $t	\geq 0,$} $B(t)= B_t $ denotes the standard Brownian
motion (BM), \ep{such that $B(0)=0$ with probability one (w.p.1),} then
\begin{equation} \label{gaussmarkov}
Y(t) = m(t) + h_2(t) B(r(t)), \ t \ge 0,
\end{equation}
is a  continuous GM process with mean
$m(t)$ and factorizable covariance $c(s,t)= h_1(s) h_2(t), $ for
$0 \le s \le t .$ \par\noindent
\ep{By substituting in \eqref{fractintegral} the function $f(t)$ with $Y(t)$ process,} our aim is to study the
fractional RL \ep{(pathwise)} integral of the GM process $Y(t),$ namely the process
\begin{equation}
X^ \alpha (t) = {\cal I }^ \alpha  (Y) (t) = \frac {1 } {\Gamma (\alpha) } \int _0 ^t (t-s)^{ \alpha -1} Y(s)  ds .
\end{equation}
\ep{Note that, due to the almost surely continuity of paths of the GM process $Y(t)$, the integral process $X^ \alpha (t)$ is well defined and it is adapted 
in the same probability space of $Y(t)$.}
\par
\ep{Referring to the neuronal model \eqref{1}-\eqref{2}, assuming that $V(0)=0$ (and, in some cases, also $\upeta(0)=0$), the RL fractional integral ${\cal I }^ \alpha$ 
is used as the left-inverse of the Caputo derivative ${\cal D}^\alpha$ (see, \cite{kil06},\cite{Mal15}). In this way, we find that the solution $V(t)$ of \eqref{1}-\eqref{2}  involves the RL fractional integral process of the GM process $\upeta(t)$, i.e., specifically }
\ep{
\begin{equation}
{\cal I }^\alpha ({\cal D}^\alpha V(t))={\cal I }^ \alpha \left(\frac{g_L}{C_m}V_L\right)+{\cal I }^ \alpha\left(\frac{{\upeta(t)}}{C_m}\right), \ \hbox{ with }V(0)=0.\\
\end{equation}}
\ep{The investigation of the last term in the right-hand-side of the  above equation has to be done. We proceed in this direction in the following.}

\subsection{The fractional Riemann-Liouville integral of BM (FIBM)}

\ep{As a first case, we consider the icon of the non-stationary GM processes: the standard Brownian motion.}
Let us consider the fractional Riemann-Liouville integral of
$B_t $ (FIBM),  that is 
$$ {\cal I }^ \alpha  (B) (t) = \frac {1 } {\Gamma
(\alpha) } \int _0 ^t (t-s)^{ \alpha -1} B(s)  ds .$$  
\ep{It has not an immediate application in the neuronal model \eqref{1}-\eqref{2}, but it will play a rule in the construction of fractional integrals of GM processes. Nevertheless, a possible model that can be suitable designed is that composed by the following coupled equations:}
\ep{
\begin{eqnarray}\label{11}
{\cal D}^\alpha V(t)&=&\frac{g_L}{C_m}V_L+\frac{{{\upeta}(t)}}{C_m}, \ V(0)=0,\\
d{\upeta(t)}&=&\frac{\varsigma}{\tau}dB(t), \hspace*{1.3cm} \upeta(0)=0, \label{22}
\end{eqnarray}}
\par\noindent
\ep{that, for $V_L=0,C_m=1$, are solved by the fractional integral process of a Brownian motion $B(t)$, i.e. 
$${\cal I }^\alpha ({\cal D}^\alpha V)={\cal I }^ \alpha{(\upeta)}=\frac{\varsigma}{\tau}{\cal I }^ \alpha{(B)}.$$}
\ep{Note that the SDE \eqref{22} is a prototype of integrate-and-fire neuronal models (\cite{Lan08}). The above model describes the neuronal membrane voltage as the fractional integral of a Brownian input process.} 
  Abut such kind of processes, the
following holds.
\begin{Proposition} \label{prop}
For $0 \le \alpha \le 1, $ let be ${\cal I }^ \alpha  (B) (t)$ the
fractional RL integral of $B_t$ of order $\alpha ;$
then, for $0 \le u \le t ,$ the covariance of ${\cal I }^ \alpha
(B) (u)$ and ${\cal I }^ \alpha  (B) (t)$ is:
\begin{equation} \label{covarianceofI(t)}
cov({\cal I }^ \alpha  (B)
(u), {\cal I }^ \alpha  (B)
(t)) = \frac 1 { \Gamma ^2 (\alpha) } \left [ \frac { t^{\alpha +1} u^ \alpha } { \alpha ^2 ( \alpha +1) }  - \frac {t H_ \alpha (u,t)}
{\alpha ( \alpha +1) } +
\frac { J_ \alpha (u,t) } {\alpha (\alpha +1)}   \right ]
\end{equation}
where
\begin{equation} \label{definitionofJ}
J _ \alpha (u,t) = \int _0 ^u s (u-s) ^{ \alpha -1} (t-s) ^\alpha ds ,
\end{equation}
\begin{equation} \label{definitionofH}
H _ \alpha (u,t) = \int _0 ^u  (u-s) ^{ \alpha -1} (t-s) ^\alpha ds .
\end{equation}

Taking $u=t$ in \eqref{covarianceofI(t)}, one gets the variance of
${\cal I }^ \alpha  (B) (t):$
\begin{equation} \label{varianzaFIBM}
Var ( {\cal I }^ \alpha  (B)
(t)) = \frac 1 { \Gamma ^2 ( \alpha )} \cdot \frac {t^{ 2 \alpha +1} } {\alpha ^2 ( 2 \alpha +1 ) } = \frac {t^{2 \alpha +1} } {(2 \alpha +1) \Gamma ^2 ( \alpha +1) }.
\end{equation}
\end{Proposition} \label{Propuno}
{\it Proof.}
Since $E(B(s) B(v))= \min (s,v),$  we have for $u \le t:$
$$
 cov({\cal I }^ \alpha  (B) (u), {\cal I }^ \alpha  (B)(t))= \frac 1 { \Gamma ^2 ( \alpha)} E \left ( \int _0 ^u (u-s)^ {\alpha -1} B(s) ds \cdot  \int _0 ^t  (t-v)^ {\alpha -1}B(v) dv \right )
$$
\begin{equation} \label{formulacov}
= \frac 1 { \Gamma ^2 ( \alpha)} \int _0 ^u ds \int _0 ^t dv (u-s)^ {\alpha -1} (t-v)^ {\alpha -1} \min (s,v);
\end{equation}
the \ep{last} integral can be split into three parts:
$$
 \int _0 ^u ds \int _0 ^s dv (u-s)^ {\alpha -1} (t-v)^ {\alpha -1} v + \int _0 ^u ds \int _s^u dv (u-s)^ {\alpha -1} (t-v)^ {\alpha -1} s
 $$
\begin{equation} \label{sommaintegrali}
+  \int _0 ^u ds \int _u^t dv (u-s)^ {\alpha -1} (t-v)^ {\alpha -1} s = I_1 + I_2 + I_3,
 \end{equation}
 where, for simplicity of notations we let drop the dependence of $I_i$ on $\alpha .$ \par\noindent
We have:
$$I_1= \int _0 ^u (u-s)^ {\alpha -1} ds \int _0 ^s  (t-v)^ {\alpha -1} v dv ;$$
the inner integral is equal to
$$ \frac {(t-s)^ {\alpha +1} } {\alpha +1 } - \frac {t (t-s)^ {\alpha} } {\alpha }
+ \frac { t^ {\alpha +1} } {\alpha ( \alpha +1)}.$$
Then, by straightforward calculation, we obtain
$$ I_1 = \frac {t^ {\alpha +1} u^ \alpha } {\alpha ^2 ( \alpha +1) }  - \frac { t} {\alpha ( \alpha +1) } \int _0^u (u-s)^ { \alpha -1} (t-s) ^ \alpha ds
- \frac 1 { \alpha +1}  \int _0 ^u s (u-s) ^{ \alpha -1} (t-s) ^\alpha ds $$
$$= \frac {t^ {\alpha +1} u^ \alpha } {\alpha ^2 ( \alpha +1) } - \frac {t } {\alpha (\alpha +1) } H_ \alpha (u,t) -
\frac {1} {\alpha +1 } J_ \alpha (u,t),$$
where $J_ \alpha (u,t)$ and $H_ \alpha (u,t)$ are given by \eqref{definitionofJ} and \eqref{definitionofH}.
As for the second integral, we have:
$$ I_2= \frac 1 \alpha \left [ \int _0 ^u s(u-s)^ {\alpha -1} (t-s) ^ \alpha ds - (t-u) ^ \alpha \int _ 0 ^u s (u-s) ^ { \alpha -1} ds \right ] $$
$$= \frac 1 \alpha J_ \alpha (u,t) - \frac 1 \alpha  (t-u) ^ \alpha \frac {u^ {\alpha +1} } { \alpha ( \alpha +1 )} .$$
The third integral is:
$$ I_3 =  \left ( \int _0 ^u s (u-s)^ {\alpha -1} ds \right ) \left ( \int _u^t  (t-v)^ {\alpha -1} dv \right ) =
\frac {u^ {\alpha +1} } {\alpha ^2 ( \alpha +1) } \cdot (t-u) ^ \alpha .$$
Finally, by summing the three integrals, and inserting in \eqref{sommaintegrali},  we obtain \eqref{covarianceofI(t)}.
\par
To obtain the variance, it is enough to put $u=t$ in \eqref{covarianceofI(t)} and to note that
$J_ \alpha (t,t) = \frac { t^ {2 \alpha +1}} {2 \alpha (2 \alpha +1) } $  and
$H_ \alpha (t,t) = \frac { t^ {2 \alpha }} {2 \alpha  } ;$ then,
 \eqref{varianzaFIBM} follows.
\par
\hfill $\Box$
\bigskip

\begin{Remark}
For $\alpha =1$ one has $J_1(u,t)= t u^2/2 - u^3/3$ and $H_1(u,t)= t^2/2 - (t-u)^2/2,$   so from  \eqref{covarianceofI(t)}
one gets the well-known  covariance function of integrated BM (see e.g \cite{ross:2010}):
\begin{equation} \label{covintegratedBM}
cov \left ( \int _0^u B_s ds , \int _ 0 ^t B_v dv \right )= u^2 \left (\frac t 2 - \frac u 6 \right );
\end{equation}
moreover, from \eqref{varianzaFIBM},
its variance follows:
\begin{equation}
var \left ( \int _0^t B_s ds \right )= t^3/3.
\end{equation}
\end{Remark}
\bigskip

\subsubsection{Numerical evaluations: the FIBM covariance}
We have computed $cov({\cal I }^ \alpha  (B) (u), {\cal I }^ \alpha  (B)(t))$ for various values of $\alpha$ and $u \le t;$ since
$J_ \alpha (u,t)$ and $H_ \alpha (u,t)$ cannot be obtained in closed form for any $\alpha ,$ their value has been obtained by numerical integration.
\par
\begin{figure}
\centering
\includegraphics[height=0.3 \textheight]{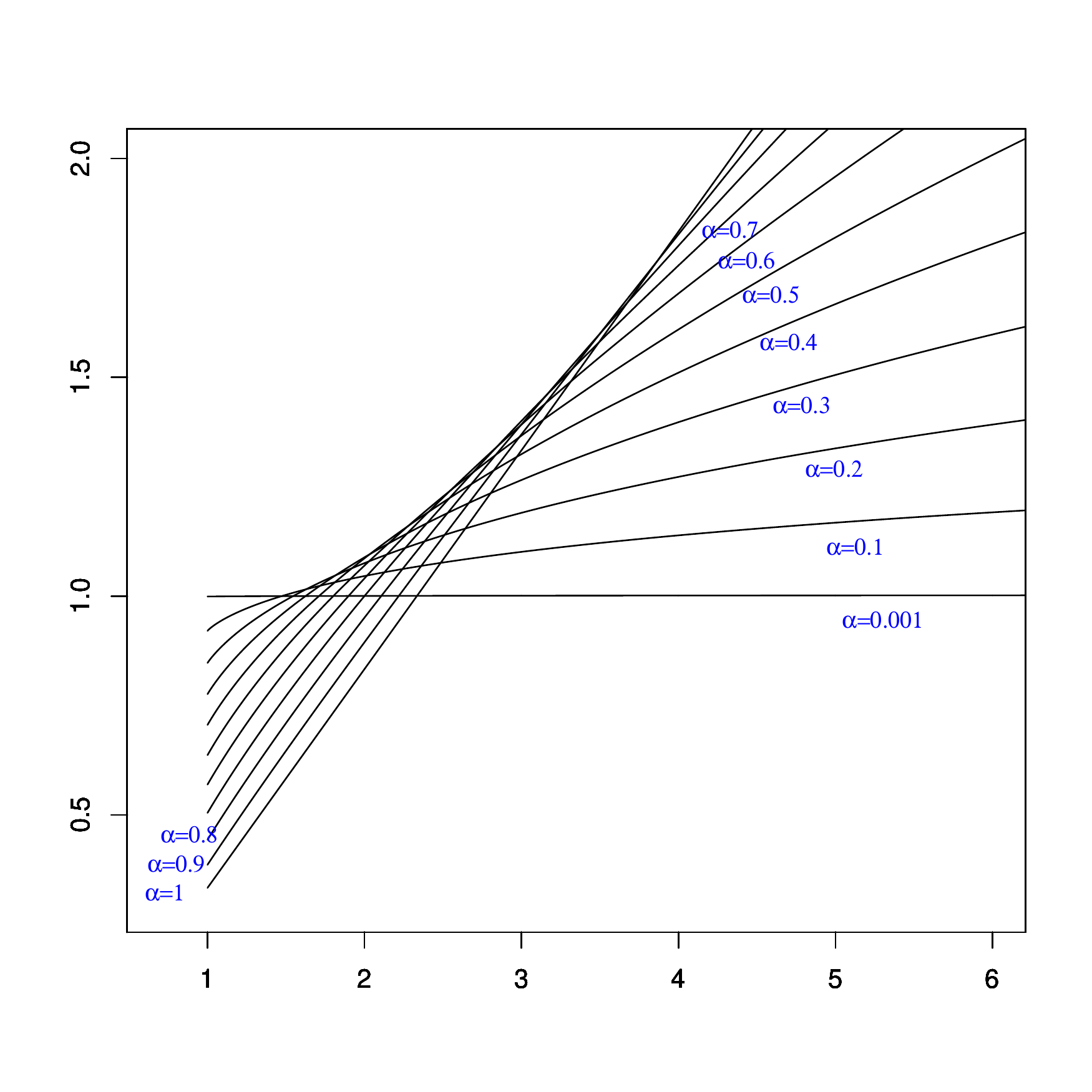}
\caption{Plot of
the covariance function of FIBM,
$cov({\cal I }^ \alpha  (B) (u), {\cal I }^ \alpha  (B)(t)),$
given by \eqref{covarianceofI(t)},
as a function of $t\ge u,$ for $u=1$ fixed, and  various values of
$\alpha \in (0,1).$
The horizontal line represents the
covariance calculated for $\alpha =0;$ the curve corresponding to
$\alpha =1$ matches the straight line of equation $y= t/2 - 1/6 .$
}
\label{covFIBM}
\end{figure}

\noindent
In the Figure \ref{covFIBM}, we report the graph of $cov({\cal I }^ \alpha  (B) (u), {\cal I }^ \alpha  (B)(t))$
 given by \eqref{covarianceofI(t)}, as a
function of $t \ge u,$ for $u=1$ fixed, and  various values of $\alpha
\in (0,1).$
\par
\begin{figure}
\centering
\includegraphics[height=0.3 \textheight]{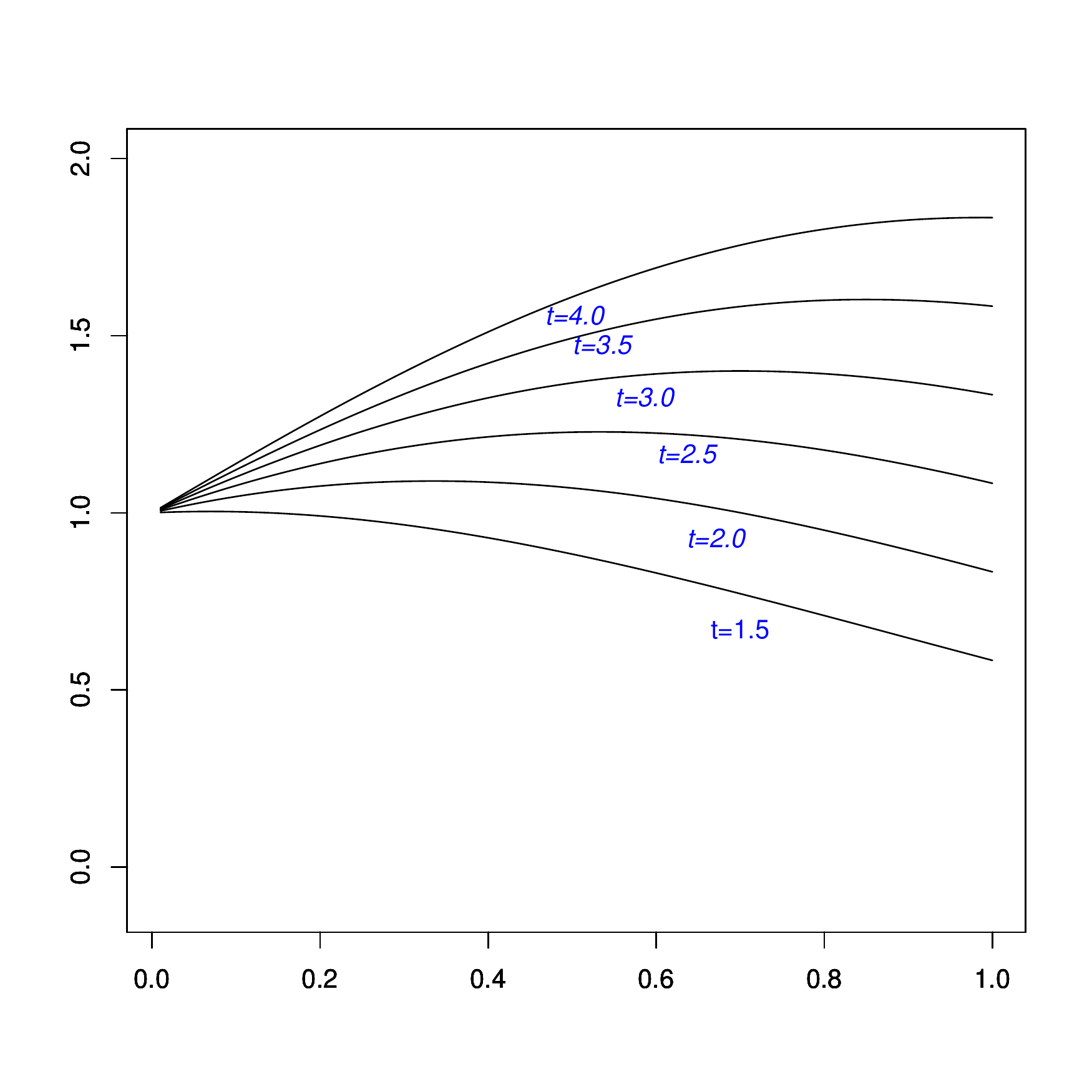}
\includegraphics[height=0.3 \textheight]{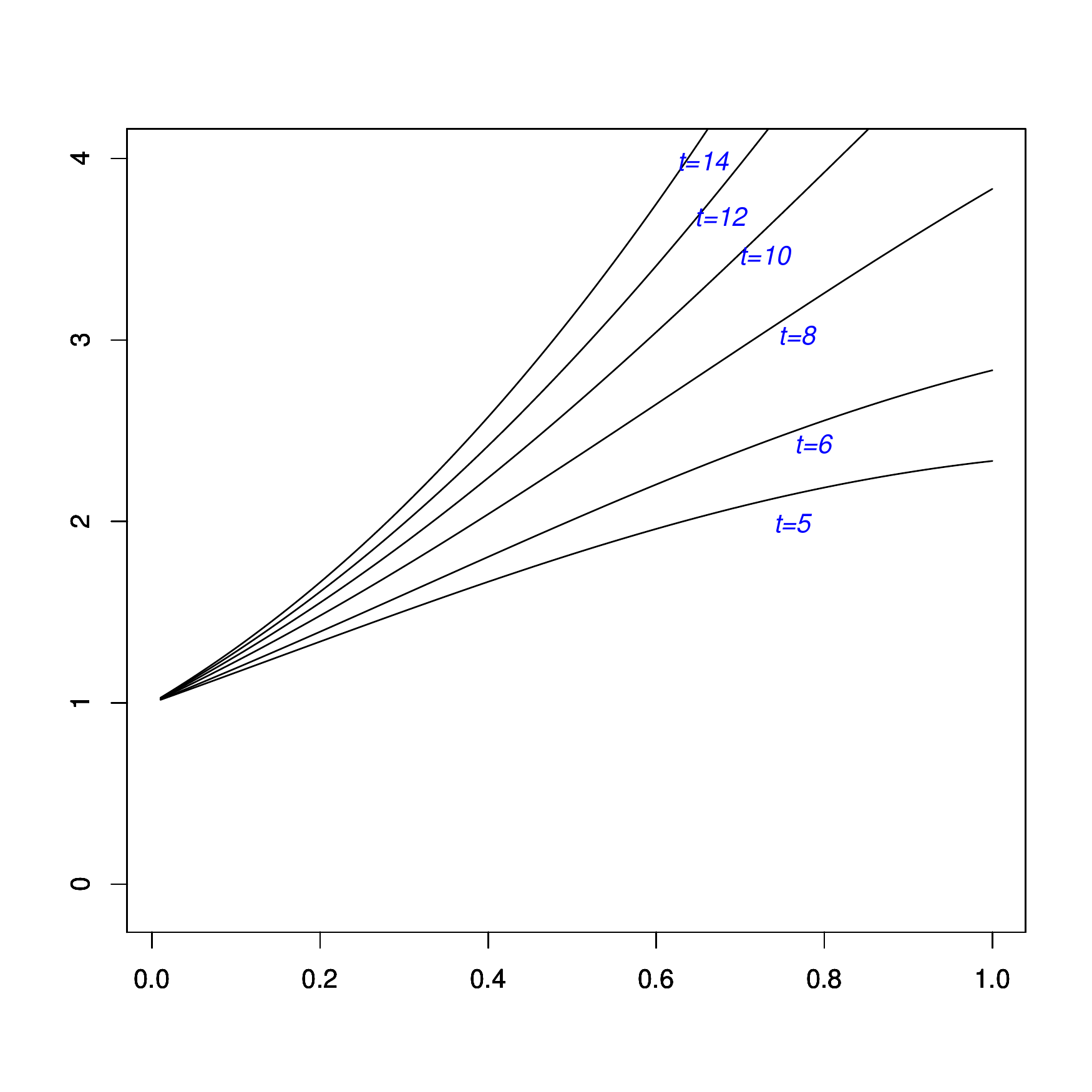}
\caption{
Plot of the covariance function of FIBM,
given by \eqref{covarianceofI(t)}, as a function of $\alpha \in
(0,1),$ for $u=1$ and  $ t = 1.5, \ 2, \ 2.5, \ 3, \ 3.5, \ 4 $ (left panel); $ t = 5, \ 6, \ 8, \ 10, \ 12, \ 14 $ (right panel). The larger $t,$ the higher the curve.  }
\label{covFIBMbis}
\end{figure}

\noindent
In the Figures \ref{covFIBMbis} we report the graph of the same
covariance as a function of $\alpha \in (0,1),$ for $u=1$ and various values of $t \ge 1.$
\par\noindent
It appears that, for fixed $u$ and $\alpha$  the covariance function of ${\cal I }^ \alpha  (B) (t)$
increases, as a function of $t \ge u .$ For fixed $u,$ there exists a value $\bar t _u $ such that for $ t < \bar t _u$  the covariance
function is decreasing as a function of $\alpha ,$  while for $t > \bar t _u $ it is increasing  (for $u=1$ one has
$ \bar t _u \simeq 1.5 ,$ see Figure \ref{covFIBM}). 
\par\noindent
Moreover,
for fixed $u$ and $t,$   a value $\bar \alpha _t $ exists at which the covariance attains its maximum, that is, for $\alpha \le \bar \alpha _t $ the
covariance function is increasing, as a function of $\alpha,$ while it is decreasing for $\alpha > \bar \alpha _t$ (see Figure \ref{covFIBMbis}).
 \par\noindent
Of course, since ${\cal I }^ 0  (B)(t) = B(t),$ for $1=u \le t$ one has 
$$cov({\cal I }^ 0  (B) (u), {\cal I }^ 0 (B)(t))= cov (B(u), B(t)) = \min (u,t) = 1 .$$
Notice that the value $\alpha =0$ cannot be directly substituted in
\eqref{covarianceofI(t)}, because $ J_ \alpha (u,t), \  H_ \alpha (u,t)$ and $ \Gamma ( \alpha)$ diverge, as $ \alpha \rightarrow 0 ^+.$
Thus, numerical calculation allows only to evaluate the covariance for small, but positive $\alpha$ (in fact, we have taken $\alpha = 0.001).$
\par\noindent
For $\alpha =1,$ we recover the covariance of integrated BM; in fact, for $u=1$ it is equal to $t/2 - 1/6 $
(see \eqref{covintegratedBM})
whose graph matches the curve obtained for $\alpha =1 .$ \par

\section{The fractional Riemann-Liouville integral of a Gauss-Markov process}
Let us consider the GM process $Y(t)$ given by
\eqref{gaussmarkov}; we have
$$
{\cal I} ^\alpha ( Y) (t) = \frac 1 {\Gamma ( \alpha )} \int _ 0 ^t (t-s) ^ { \alpha -1} Y(s) ds =
\frac 1 {\Gamma ( \alpha )} \int _ 0 ^t (t-s) ^ { \alpha -1} [ m(s) + h_2 (s) B ( r(s)) ] ds
$$
\begin{equation} \label{fractalintegralGM}
= \frac 1 {\Gamma ( \alpha )} \int _0 ^t [ \widetilde m (s,t) + \widetilde h_2 (s,t) B(r(s)) ] ds,
\end{equation}
where
\begin{equation} \label{mtilde}
 \widetilde m(s,t) = (t-s) ^ {\alpha -1} m(s), \ \widetilde h _2 (s,t) = (t-s) ^ {\alpha -1} h_2 (s) 
 \end{equation}
(here, the dependence on $\alpha$ was omitted to avoid an heavy notation).
Then, by straightforward calculations, \ep{and applying the linearity in $L^1$ of the fractional RL integral operator,} we obtain:
\begin{Proposition} \label{proposizionedue}
The fractional RL integral of $Y(t),$
\begin{equation}
{\cal I }^ \alpha (Y) (t) = \frac {1 } {\Gamma (\alpha) } \int _0 ^t (t-s)^{ \alpha -1} Y(s) ds ,
\end{equation}
is normally distributed with mean $M(t) = \frac 1 { \Gamma ( \alpha )} \int _0 ^t \widetilde m (s,t) ds $ and
covariance, for $0 \le u \le t:$
$$
cov \left ({\cal I }^ \alpha (Y) (u), {\cal I }^ \alpha (Y)(t) \right ) = \frac 1 {\Gamma ^2 ( \alpha)} cov \left ( \int _0 ^u \widetilde h_2 (s,t) B( r(s)) ds, \int _0 ^t \widetilde h_2 (v,t) B( r(v)) dv \right )$$
\begin{equation} \label{covFIGM}
= \frac 1 {\Gamma ^2 ( \alpha)} \int _0 ^ u ds \int _0 ^ t dv  \ \widetilde h_2 (s,t) \widetilde h_2 (v,t)  \min (r(s),r(v)).
\end{equation}

\end{Proposition}

\begin{Remark}
The integral in \eqref{covFIGM} can be calculated, by splitting it into three parts, as in the proof of Proposition \ref{prop}; thus, one obtains:
$$
cov \left ({\cal I }^ \alpha (Y) (u), {\cal I }^ \alpha (Y)(t) \right ) = \frac 1 {\Gamma ^2 ( \alpha)} \Big [ \int _0 ^u ds (u-s) ^ { \alpha -1} h_2 (s)
\int _0 ^s dv (t-v) ^ {\alpha -1} r(v) h_2 (v)
$$
$$ + \int _ 0 ^u ds (u-s)^ {\alpha -1} h_2 (s) r(s) \int _ s ^u dv (t-v) ^ { \alpha -1} h_2 (v) $$
\begin{equation} \label{covFIGMespl}
+ \left ( \int _ 0 ^u ds (u-s)^ {\alpha -1} r(s) h_2 (s) \right ) \left ( \int _u^ t  dv (t-v) ^ {\alpha -1} h_2 (v) \right ) \Big ] :=
\frac 1 {\Gamma ^2 ( \alpha)} \left [ \widetilde I_1 + \widetilde I_2 + \widetilde I_3 \right ],
\end{equation}
where, for simplicity of notations we let drop the dependence of $\widetilde I _i$ on $\alpha .$
For $u=t,$ $\widetilde I_3$ vanishes and $\widetilde I_2 = \widetilde I_1,$ so one has:
\begin{equation} \label{varFIGMespl}
Var \left ({\cal I }^ \alpha (Y) (t) \right ) =
\frac 2 {\Gamma ^2 ( \alpha)} \int _0 ^t ds (t-s) ^ { \alpha -1} h_2(s) \int _0 ^s dv (t-v)^ { \alpha -1} r(v) h_2(v).
\end{equation}

\end{Remark}
\section{Fractional integral of GM processes and comparisons with the ordinary integral of  GM processes}
In this section, we study the qualitative behavior of the fractional integral over time  of a GM process, when varying the
order $\alpha  \in (0,1),$  and we compare it with the corresponding ordinary integral. Computer simulation is also used to illustrate these behaviors.
First, we recall from \cite{abupir:physicaA17} the following result, regarding the ordinary integral over time of a GM process.

\begin{Theorem}   \label{proposition1}
Let $Y$ be a GM process of the form \eqref{gaussmarkov} with $r (0)=r_0
\ge 0 ,$
$r'(t)>0$ for $t >0;$
then $X(t)=x + \int _ 0 ^t Y(s) ds $ can be written as $\mathcal V(t) + \eta \int _0 ^t h_2
(s) ds,$ in which $\eta = B(r(0))$ and $\mathcal V(t)$
is normally distributed with  mean $x+ M(t)$  and variance $\gamma (\rho(t)),$ where
$\rho (t)= r(t)- r(0),$
$M(t)= \int _0 ^t m(s) ds,$
$\gamma (t)= \int _0 ^t (R(t) - R(s)) ^2 ds $ and $R (t)= \int _0 ^t h_2(\rho ^{-1}
(s))/ \rho '( \rho ^{-1} (s)) ds .$

\end{Theorem}

\hfill $\Box$
\bigskip

 \ep{Note that, generally the $\mathcal  V(t)$ process and the r.v. $\eta$ are not independent, when $r(0)\ne 0$ (if $r(0)=0,$ then the term $\eta$ vanishes).}

\subsection{Fractional integrated BM (FIBM)}
If $Y(t)= B_t,$ taking $m(t)=0, \ h_2(t)=1, \ r(t) = \rho (t)= t$ in Theorem \ref {proposition1}, we have for the ordinary integral of BM:
\begin{equation}
\int _0 ^ t B_s ds \sim {\cal N } (0, t^3/3) ,
\end{equation}
where ${\cal N } (\mu, \sigma ^2) $ denotes a Gaussian  r.v. with
mean $\mu$ and variance $ \sigma ^2 .$  \par
By using Proposition
\ref{prop}, we obtain that the fractional integral of order $\alpha,$
i.e. 

$${\cal I} ^ \alpha  (B) (t)= \frac 1 { \Gamma ( \alpha)} \int
_0 ^t (t-s)^ { \alpha -1} B_s ds \sim {\cal N}
(0, t^ {2 \alpha +1} /[\alpha ^2 ( 2 \alpha +1)\Gamma ^2 ( \alpha
)]).$$

\par\noindent Thus, we have: \vspace{0.3cm}
\par\noindent $\bullet$ for $\alpha \rightarrow
0^+, \ {\cal I} ^ \alpha  (B)(t) \sim
{\cal N} (0, t ),$ because $\alpha \Gamma ( \alpha) = \Gamma ( \alpha +1)$
and so $\alpha ^2 ( 2 \alpha +1) \Gamma ^2 ( \alpha) $ tends to $\Gamma ^2 (1) =1;$ therefore,
 ${\cal I} ^ 0  (B)(t)$ has the same distribution as $B_t.$
\vspace{0.3cm}
 \par\noindent
$\bullet$ for $\alpha =1/2,$ one has ${\cal I} ^ \alpha  (B)(t) \sim
{\cal N} (0, 2t^2/ \pi ),$ because $\Gamma ( 1/2) = \sqrt \pi .$
\vspace{0.3cm}
\par\noindent $\bullet$ for $\alpha
\rightarrow 1^- ,$ one has 
${\cal I} ^ \alpha (B) (t) \sim {\cal N}
(0, t^3/3),$ because the fractional integral coincides with the
ordinary integral.
\par\noindent

\subsubsection{Numerical evaluations of the FIBM variance and simulation paths}
\begin{figure}
\centering
\includegraphics[height=0.3 \textheight]{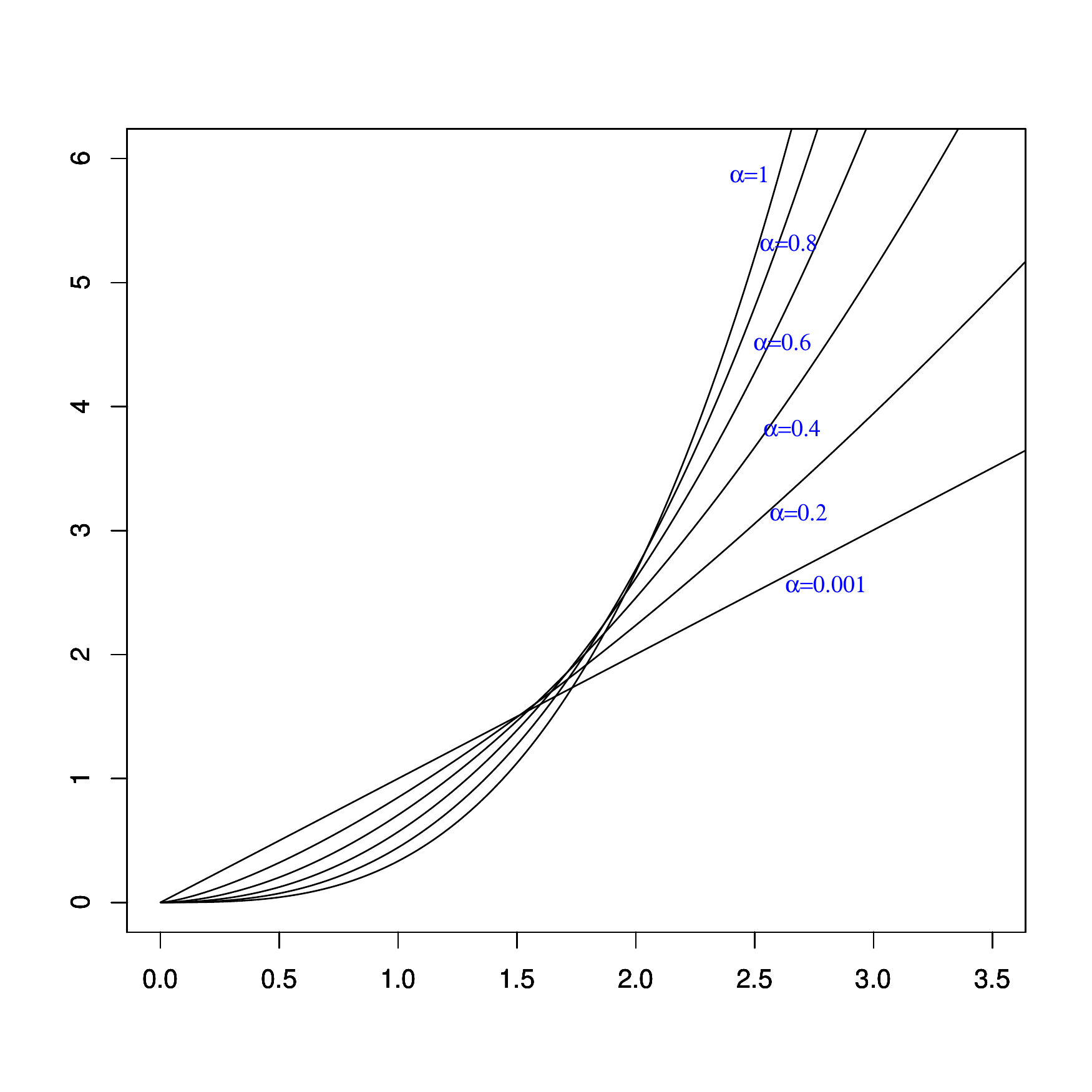}
\caption{Plot of the variance of FIBM, ${\cal I} ^ \alpha (B) (t),$  given by \eqref{varianzaFIBM}, as a function of $t \in [0,3.5],$
for various values of $\alpha \in (0,1).$ The time zone at which an inversion of behavior is observed goes from
$t\simeq  1.73 $ to $t \simeq 1.77.$ For $\alpha =0.001$ the curve is close to the graph of the function $v(t)=t,$ since for $\alpha =0$ FIBM
becomes BM; for $\alpha =1$ the curve matches the graph of the function $v(t)=t^3/3,$ since  FIBM
becomes the ordinary integral of BM.
}
\label{varFIBM}
\end{figure}

In the Figure \ref{varFIBM}  we report the shape of the variance of ${\cal I} ^
\alpha (B) (t),$ i.e. $t^ {2 \alpha +1} /[( 2 \alpha
+1)\Gamma ^2 ( \alpha +1)]$, as a
function of $t>0,$ for various values of $\alpha \in (0,1).$ As we
see, for fixed $\alpha$ the variance is increasing, as a function
of $t.$ Moreover, for small enough values of  $t,$ the curves become ever lower, as $\alpha$ increases,
that is, the variance decreases as a function of $\alpha;$ for  larger enough values of $t$  this behavior is overturned, because
the variance increases with $\alpha.$ The time zone at which an inversion of behavior is observed goes from
$t\simeq  1.73 $ to $t \simeq 1.77.$ For $\alpha =0.001$ the curve is close to the graph of the function $v(t)=t,$ since
for $\alpha =0 $ FIBM
becomes BM; for $\alpha =1$ the curve matches the graph of the function $v(t)=t^3/3,$ since  FIBM
becomes the ordinary integral of BM.

\begin{figure}
\centering
\includegraphics[height=0.3 \textheight]{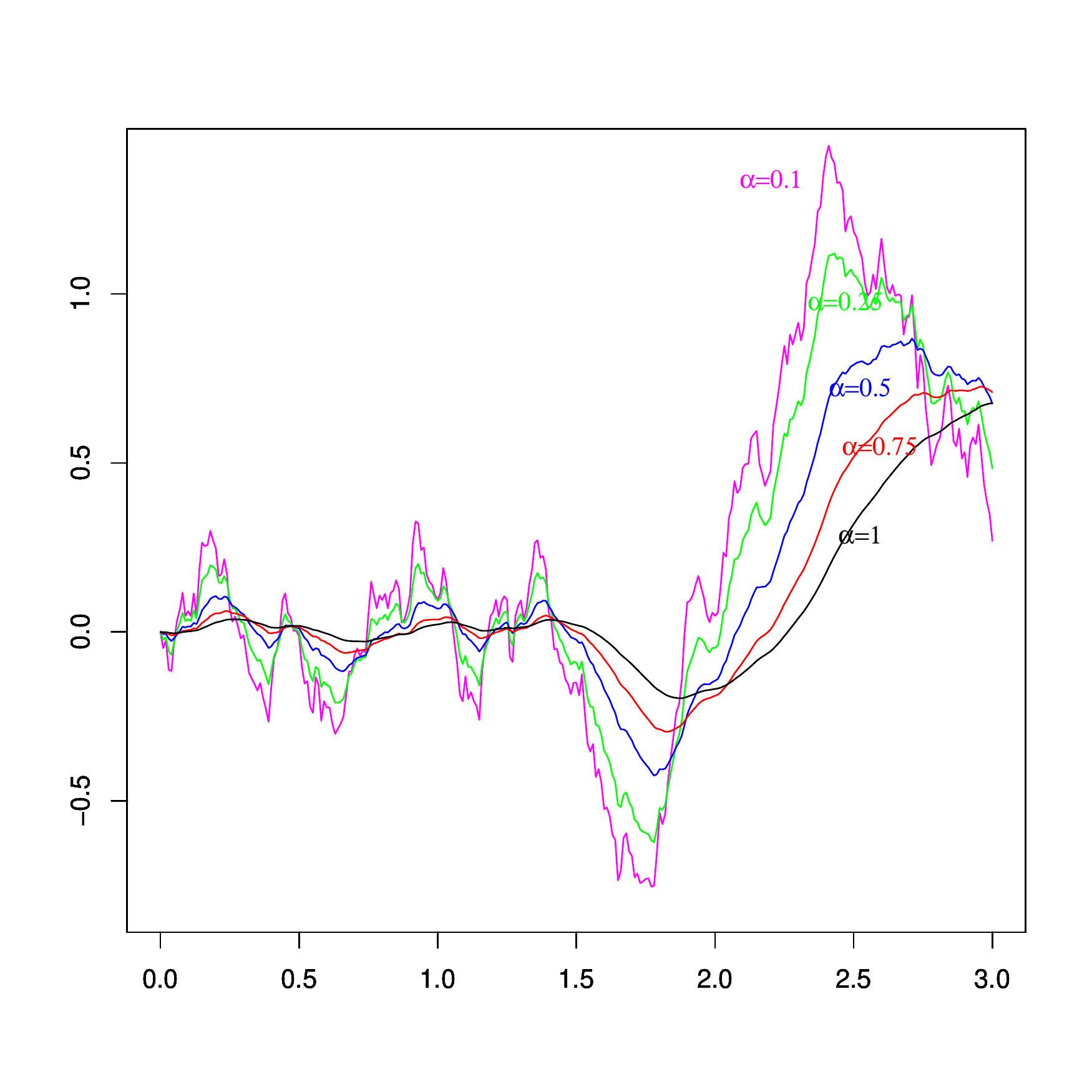}
\caption{Simulated trajectories of FIBM as function of time $t,$
for some values of $\alpha \in (0,1).$
 }
 \label{simuFIBM}
\end{figure}

\noindent
In the Figure \ref{simuFIBM}
we report the graphs of simulated
trajectories of FIBM as function of time $t,$
for some values of $\alpha \in (0,1).$ 

\begin{figure}
\begin{minipage}[h]{1\textwidth}
\centering
{\includegraphics[width=0.33\textwidth]{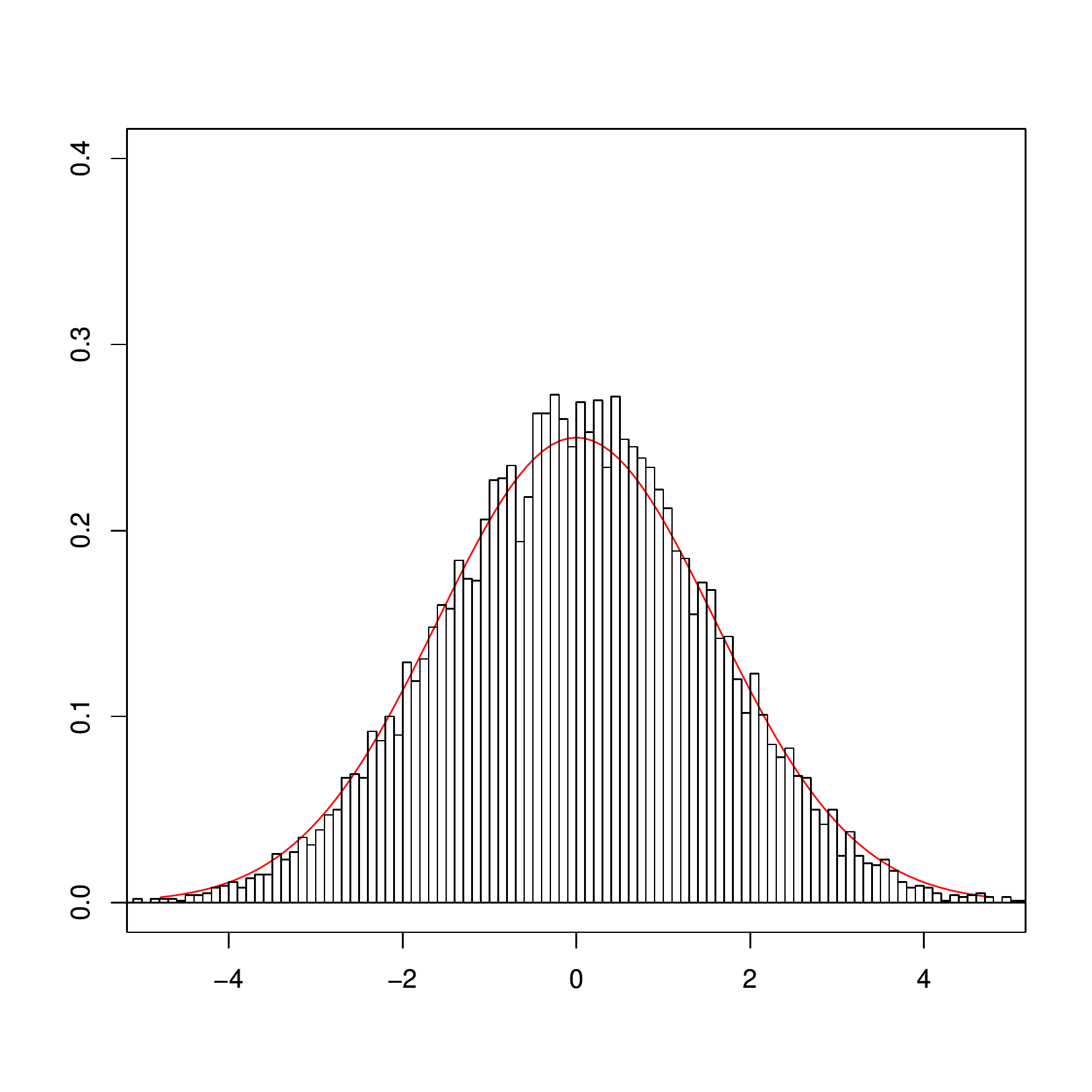}}
{\includegraphics[width=0.33\textwidth]{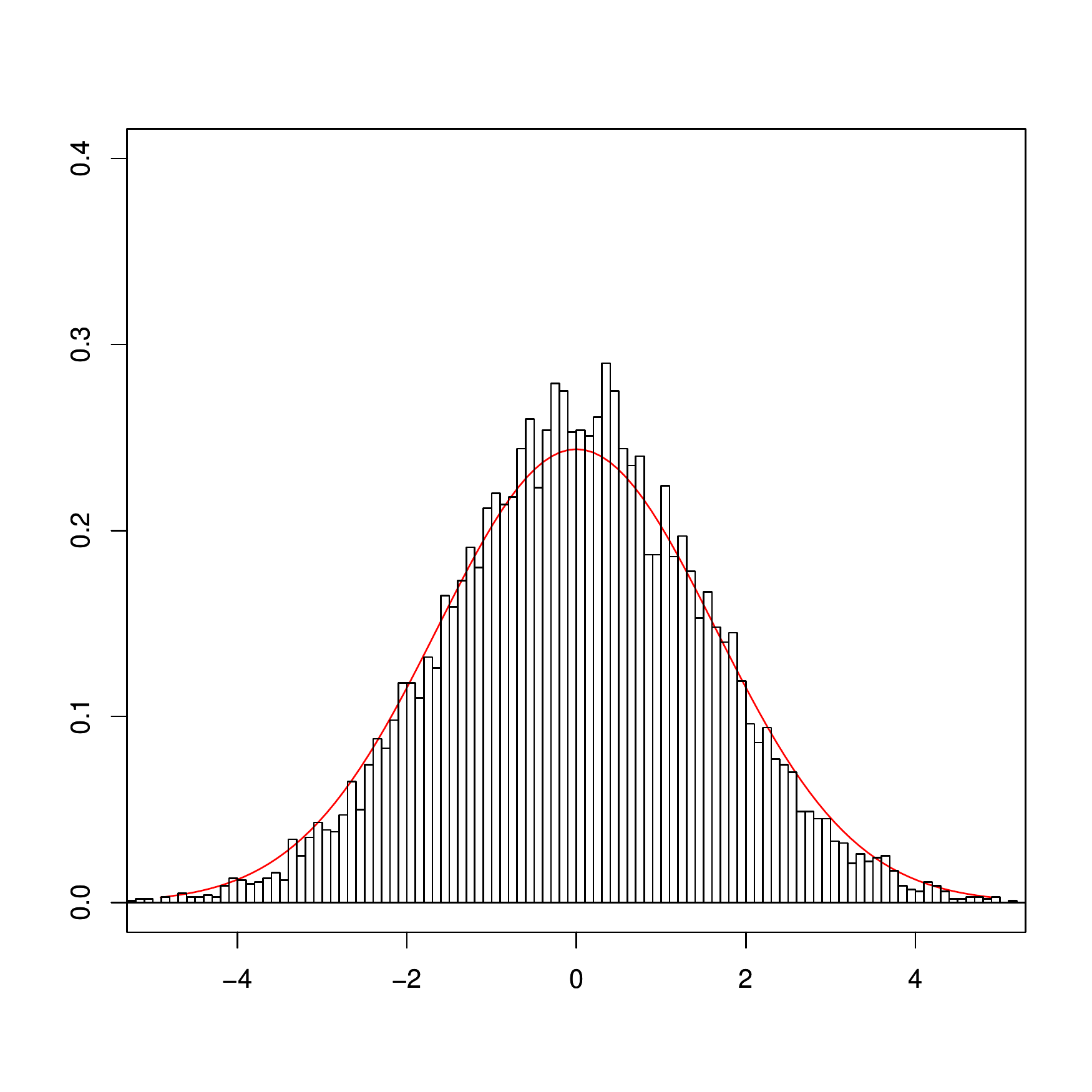}}\\
{\includegraphics[width=0.33\textwidth]{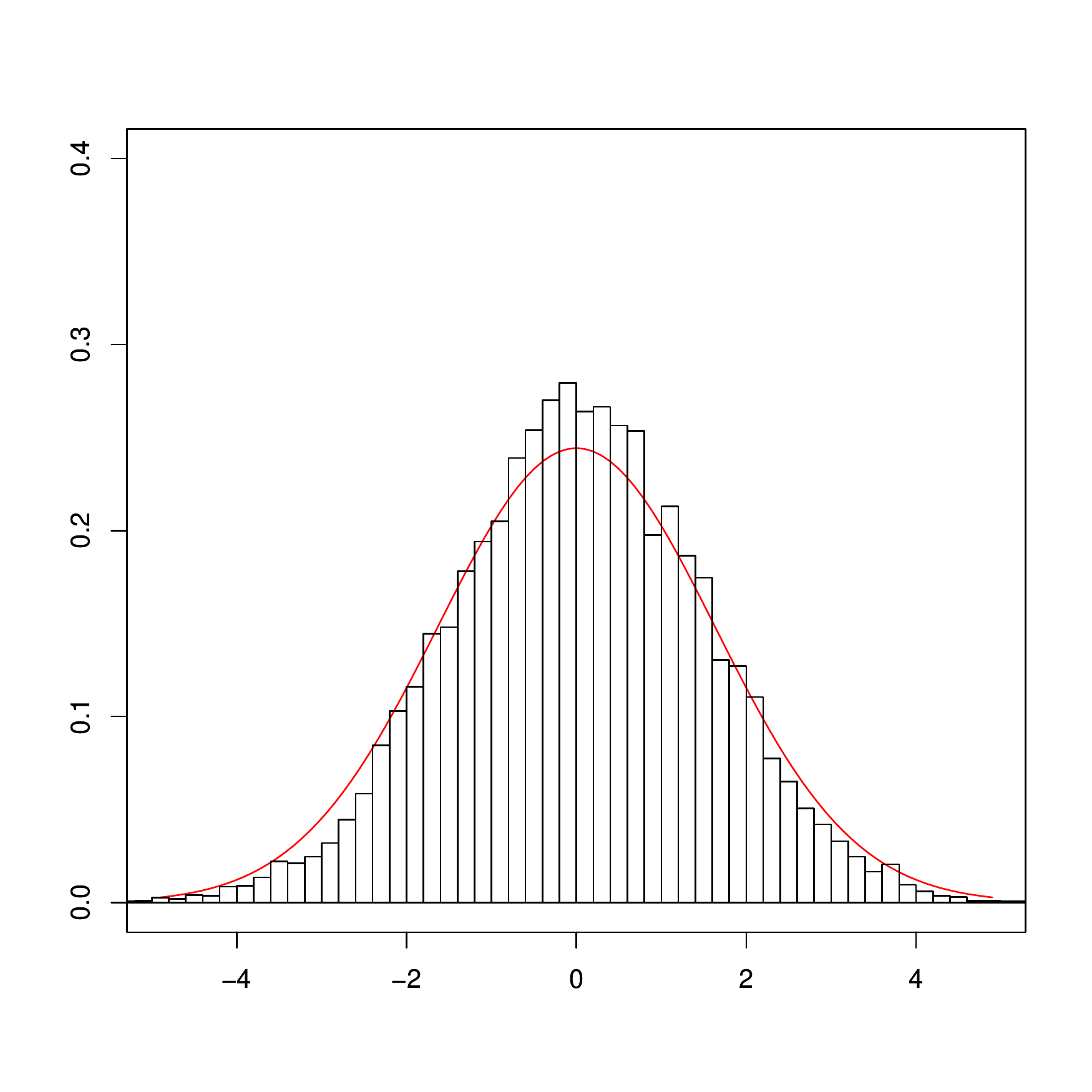}}
{\includegraphics[width=0.33\textwidth]{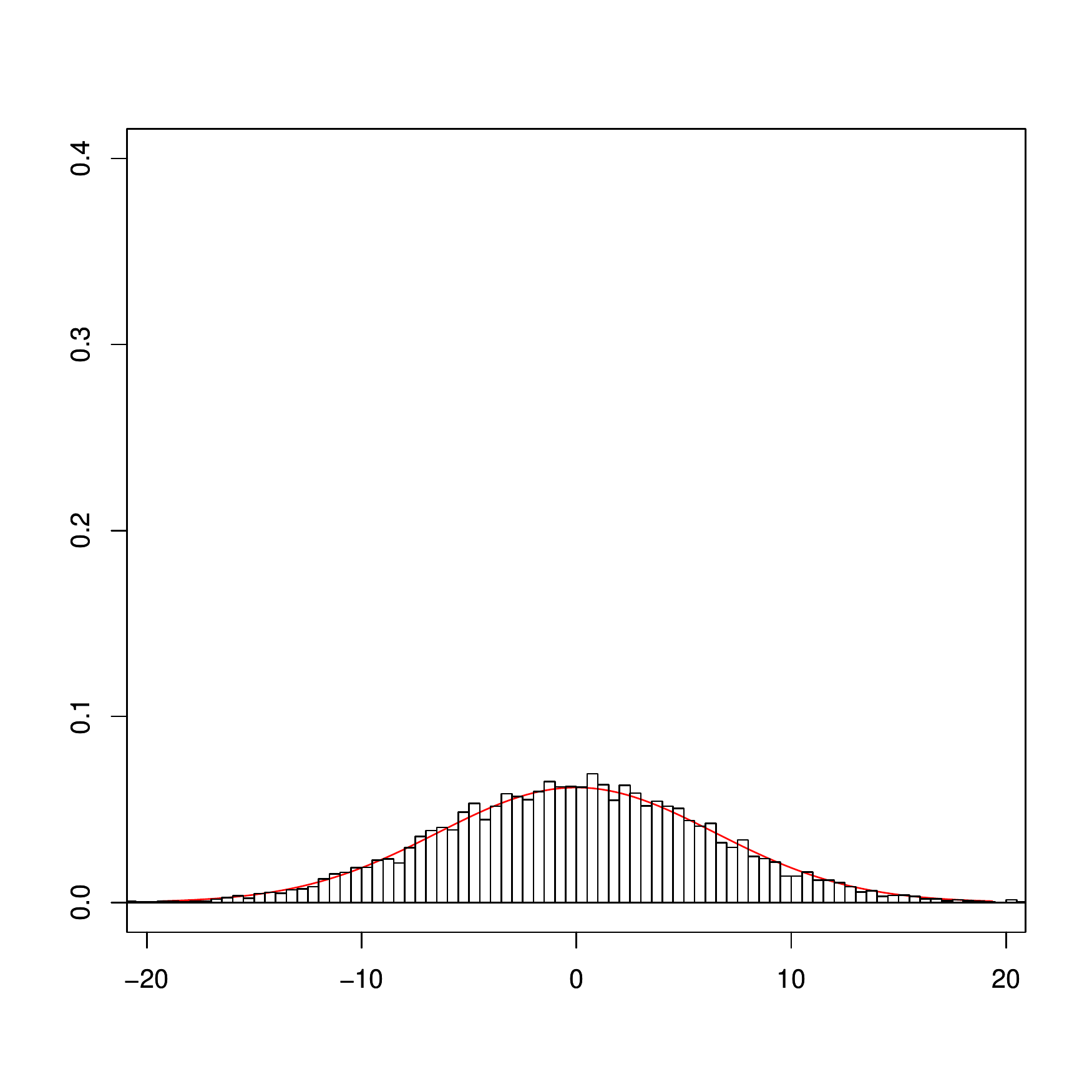}}
\caption{On top: Comparison between the density histogram of simulated values
of FIBM at $t=2,$ and the Gaussian density with zero mean and variance given by \eqref{varianzaFIBM},
for  $\alpha =0.5$ (left) and $\alpha =0.75$ (right).
On Bottom: Comparison between the density histogram of simulated values
of FIBM with $\alpha =1,$ at $t=2,$ and the Gaussian density with zero mean and variance given by \eqref{varianzaFIBM},
(left); the same for $\alpha =1$ and $t=5$ (right).}
\label{comphistFIBM}
\end{minipage}
\end{figure}

\par\noindent
Specifically, the sample paths have been obtained by using the {\bf R} software,
with time discretization step $h=0.01$ \ep{and by means of the same sequence of pseudo-random Gaussian numbers.}
\ep{
Our simulation algorithm has been realized as an R script. Referring to algorithms for the generation of pseudo-random numbers (see, for instance, \cite{Hau17}), the main steps of implementation are the following: 
\begin{itemize}
  \item calculation of the elements of $N \times N$ covariance matrix $C(t_i,t_j)$ at times $t_i,\ i=1,\ldots,N,$ of an equi-spaced temporal grid;
  \item application of the Cholesky decomposition algorithm to the covariance matrix $C$ to obtain a lower triangular matrix $L(i,j)$ such that $C=LL^T;$
  \item generation of a N-dimensional array $\bold z$ of standard pseudo-Gaussian numbers;
  \item construction of the sequence of simulated values of the correlated fractional integrated process as the array $\bold{x}=L\bold{z}.$  
   \end{itemize}
}
\par\noindent

Finally, the array $\bold x$ provides the simulated path $(X(t_1),\ldots,X(t_N)),$ which components have the assigned covariance.
\par\noindent
In particular, in the Figure \ref{comphistFIBM}, we compare the density histogram of simulated values
of FIBM at a certain time $t,$ with the Gaussian density with zero mean and variance given by \eqref{varianzaFIBM},
for various values of $\alpha .$ As we see, the matching is very good.

%
%
%

\subsection{Fractional integrated OU (FIOU) process}
Let be $Y_{OU}(t)$ the (non stationary) OU process, given by
the solution of the SDE:
\begin{equation}\label{sdeYou}
dY_{OU}(t)= - \mu (Y_{OU}(t)-\beta)dt + \sigma dB_t, \ Y_{OU}(0)=y ,
\end{equation}
where $\mu , \sigma >0$ and $\beta \in (- \infty, + \infty ).$
The explicit solution is (see e.g.  \cite{abundo:stapro12}):
\begin{equation} \label{espicitformulaOU}
Y_{OU}(t)=  \beta + e^{ - \mu t } [ y-\beta +  \widetilde B(\rho (t))] ,
\end{equation}
where $\widetilde B (t)$ is standard BM and
\begin{equation} \label{rhoforOU}
r(t)= \rho (t)= \frac {\sigma ^2} {2 \mu} \left (e ^{2 \mu
t } -1 \right )
\end{equation}
(notice that $r(0)=0).$ So, $Y_{OU}(t)$ is a GM process with:
\begin{equation} \label{mforIOU}
m(t)=  \beta + e^{ - \mu t } (y- \beta ),
\end{equation}
\begin{equation} \label{hforIOU}
h_1(t)= \frac {\sigma ^2} {2 \mu} \left (e ^{\mu
t } -e^{- \mu t } \right ), \
h_2(t)= e^{ - \mu t }
\end{equation}
and covariance
\begin{equation} \label{covOU}
c(s,t)= h_1(s) h_2(t)=
\frac {\sigma ^2} {2 \mu} \left (e^ {- \mu (t-s)} - e^ { -\mu (s+t)} \right ), \ 0 \le s \le t.
\end{equation}

\par\noindent
\ep{
Referring to the neuronal model, i.e. to the equations \eqref{1}-\eqref{2}, the OU process $Y_{OU}(t)$ stands for the process $\upeta(t)$ solution of \eqref{2} if  one takes $\tau=1/\mu$, $\varsigma/\tau=\sigma$ and a constant input current $\beta$ in place of the $I(t)$.}

By calculating the various quantities in Theorem \ref{proposition1}, we get (see
\cite{abundo:smj13}):
\begin{equation}
M(t)= \int _0 ^t \left ( \beta + e^{ - \mu s } (x- \beta ) \right ) \ ds= \beta t +
\frac {(x - \beta)} \mu  \left (1 - e^{ - \mu t } \right ) ,
\end{equation}
\begin{equation}
 \rho ^{-1} (s)= \frac 1 { 2 \mu} \ln \left (1 + \frac{ 2 \mu} {\sigma ^2} s \right
) , \ R(t) = \int _0 ^t  e^{- \mu \rho ^{-1} (s)} ( \rho ^{-1})'(s) ds = \frac {1-
e ^{- \mu \rho ^{-1} (t) }  } { \mu} ,
 \end{equation}
$$
 \gamma(t)= \frac 1 { \mu ^2} \int _0 ^t \left ( e^{- \mu \rho ^{-1} (t) } - e^{ -
\mu \rho ^{-1} (s)} \right ) ^2 ds =
\frac 1 { \mu ^2} \int _0 ^t \left ( \frac 1 { \sqrt { 1+ 2 \mu t / \sigma ^2 } } -
\frac 1 { \sqrt { 1+ 2 \mu s / \sigma ^2 } } \right ) ^2 ds $$
\begin{equation} \label{gammaforOU}
=  \frac {\sigma ^2 t } {\mu ^2 (\sigma ^2 + 2 \mu t ) } - \frac {2 \sigma ^2 } {\mu
^3 \sqrt { 1+ 2 \mu t / \sigma ^2 } }
\left ( \sqrt { 1+ 2 \mu t / \sigma ^2 } -1 \right ) + \frac { \sigma ^2 } {2 \mu ^3
} \ln \left ( 1+ 2 \mu t / \sigma ^2 \right ) .
\end{equation}
Thus, by Theorem \ref{proposition1}, being $\eta = B(r(0))=B(0)=0,$ we obtain that
the integrated OU (IOU) process (i.e. the ordinary integral over time of the OU process)
\begin{equation}
X_{IOU}(t)= \int _0 ^t Y_{OU}(s) ds = {\cal I} ^ 1 (Y_{OU}) (t)
\end{equation}
 is normally distributed with mean
$M(t)$
and variance $\gamma (\rho(t))  .$
By calculation, one obtains:
\begin{equation} \label{varIOU}
var( X_{IOU} (t))=
\frac { \sigma ^2} { 2 \mu ^3 } \left [1- e^{- 2 \mu t } -4(1-e^{ - \mu t }) + 2 \mu t \right ], \ t \ge 0.
\end{equation}
\bigskip

\noindent As far as the FIOU, i.e. ${\cal I} ^ \alpha (Y_{OU}) (t),$ is concerned,
by Proposition \ref{proposizionedue} we get that
it has normal distribution with mean
 $M(t) = \frac 1 { \Gamma ( \alpha )} \int _0 ^t \widetilde m (s,t) ds ,$ and covariance and variance
respectively given by \eqref{covFIGMespl} and \eqref{varFIGMespl}, where $r(t)$ is
as in  \eqref{rhoforOU}, $\widetilde m, \ \widetilde h_2 $ are
the corresponding functions obtained inserting \eqref{mforIOU},
\eqref{hforIOU} in \eqref{mtilde}.
Precisely, the integrals in \eqref{covFIGMespl} turn out to be:
$$ \widetilde I_1 = \frac { \sigma ^2 } \mu \int _0 ^u ds (u-s) ^ {\alpha -1} e^ { - \mu s} \int _ 0 ^s dv (t-v) ^ {\alpha -1} \sinh ( \mu v ) ,$$
$$ \widetilde I_2 = \frac { \sigma ^2 } \mu \int _0 ^u ds (u-s) ^ {\alpha -1} \sinh( \mu s) \int _ s ^u dv (t-v) ^ {\alpha -1} e^ { - \mu v},$$
\begin{equation} \label{integraliItilde}
\widetilde I_3 = \frac { \sigma ^2 } \mu \left ( \int _0 ^u ds (u-s) ^ {\alpha -1} \sinh( \mu s) \right )
\left ( \int _ u ^t dv (t-v) ^ {\alpha -1} e^ { - \mu v} \right )
.
\end{equation}
Notice that, for $\sigma =1 $ and $\mu \rightarrow 0$ one
obtains again the integrals $I_i, \ i=1, 2, 3$ in
\eqref{sommaintegrali}, since $Y_{OU} (t)$ approaches BM. \par\noindent Taking $\alpha =1$ in the
expressions above, by \eqref{covFIGMespl} one obtains the
covariance of the ordinary integral of OU (IOU):
$$
cov({\cal I}^1(Y_{OU})(u), {\cal I}^1(Y_{OU})(t))
$$
\begin{equation} \label{covIOU}
= \frac { \sigma ^2} { 2 \mu ^3 } \left (
2 \mu u +4 e^{- \mu u } - e^{- \mu (t-u)} - e^{- \mu (t+u)} + 2 e^{ - \mu t} -2 e^{ - \mu u} -2 \right ), \ 0 \le u \le t.
\end{equation}

\subsubsection{Numerical evaluations: The FIOU variance and covariance}
As we see, for $ \alpha $ different from $0$ and $1,$ the calculations required to
find the covariance of FIOU are far more complicated than in the case of
FIBM; in fact, $\widetilde I_i$ are double
integrals which cannot be found analytically, so we have calculated
them numerically, by using the {\bf R} software. 

\begin{figure}
\centering
\includegraphics[height=0.33 \textheight]{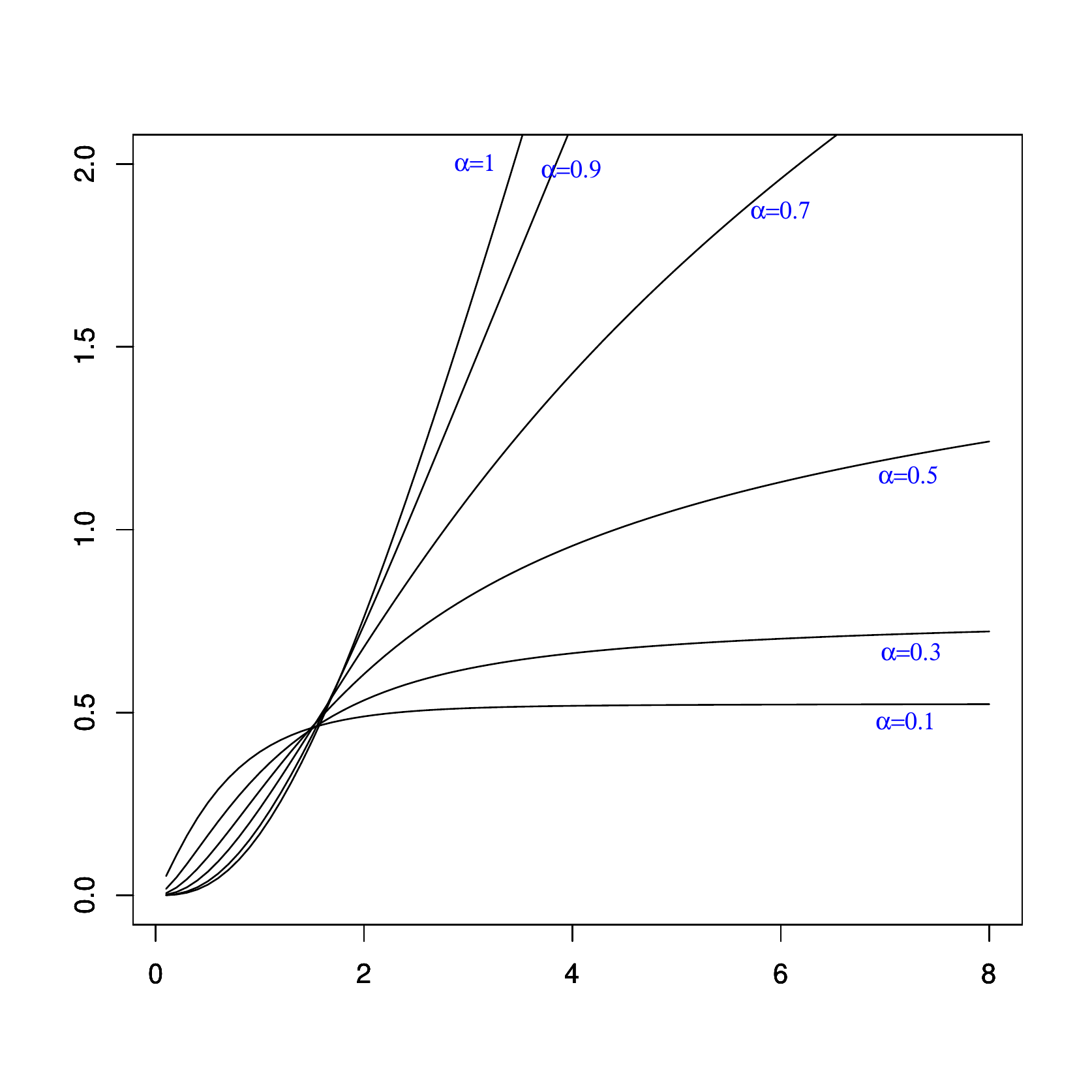}
\caption{Plot of the variance of FIOU, ${\cal I} ^ \alpha (Y_{OU}) (t)$ with $\sigma = \mu =1,$ as a function of $t \in [0,8],$
for various values of $\alpha \in (0,1).$  The time zone at which an inversion of behavior is observed goes from
$t\simeq  1.73 $ to $t \simeq 1.77.$ For $\alpha =0.1$ the curve is close to the graph of the function
$Var(t)= \frac { \sigma ^2 } { 2 \mu } (1-\exp ( - 2 \mu t ))$ (i.e. the variance of OU), since FIOU for $\alpha =0$
becomes OU; for $\alpha =1$ the curve matches the graph of the function $\gamma ( \rho (t))$
(i.e. the variance of IOU),
given by \eqref{varIOU}, since  FIOU
becomes the ordinary integral of OU.
}
\label{varFIOU}
\end{figure}

With reference to OU process with $\sigma = \mu =1$, in the Figure \ref{varFIOU}, for $u=1$ and
$u \le t \le 2,$ and several values of $\alpha,$ we
report  the shape of the variance of ${\cal I} ^
\alpha (Y_{OU}) (t).$ 
\par\noindent
As we see, for
fixed $\alpha$ the variance is increasing, as a function of $t.$
Moreover, for small enough values of  $t,$ the curves become ever
lower, as $\alpha$ increases, that is, the variance decreases as a
function of $\alpha;$ for  larger enough values of $t$  this
behavior is overturned, because the variance increases with
$\alpha.$ In fact, it appears a behavior analogous to that of the
variance of the FIBM (see Figure \ref{varFIBM}); the time zone at which an
inversion of behavior is detected, is the same observed for the variance of the
FIBM. For $\alpha =0.1$ the curve is close to the graph of the function
$v(t)= \frac { \sigma ^2 } { 2 \mu } (1-\exp ( - 2 \mu t ))$ (i.e.
the variance of OU), since for $\alpha =0$ FIOU becomes OU; for $\alpha =1$ the
curve matches the graph of the function $\gamma ( \rho (t))$ (i.e.
the variance of IOU), given by \eqref{varIOU}, since  FIOU becomes the ordinary integral of
OU. 

\begin{figure}
\centering
{\includegraphics[width=0.4\textwidth]{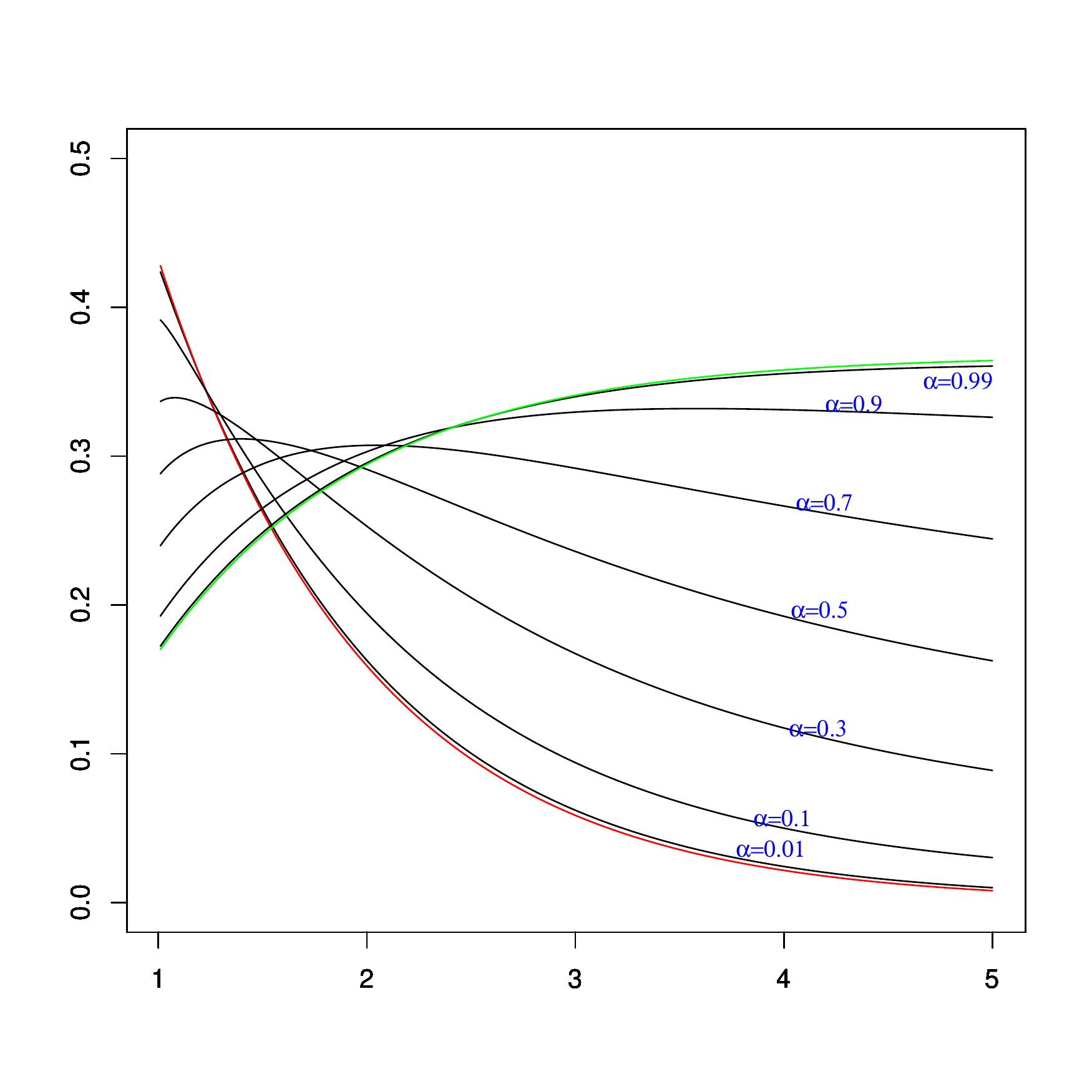}}
\caption{ Plot of the covariance of FIOU, $cov({\cal I}^\alpha(Y_{OU})(u), {\cal I}^\alpha(Y_{OU})(t)),$ with $\sigma = \mu =1,$  for  $1=u \le t \le 5,$
and various values of $\alpha \in (0,1).$
For $\alpha = 0.01$ the curve is close to the covariance of OU, given by
\eqref{covOU} (red curve), for $\alpha =0.99$ the curve is close to that of the covariance of the ordinary integral of OU (IOU),
given by \eqref{covIOU} (green curve).}
\label{covFIOU}
\end{figure}

In the Figure \ref{covFIOU} the $cov({\cal
I}^\alpha(Y_{OU})(u), {\cal I}^\alpha(Y_{OU})(t)),$ for $u=1$ and
$u \le t \le 2,$ for various values of $\alpha.$ 
We note that
for $\alpha = 0.01$ the curve is close to the covariance of OU, given by
\eqref{covOU} (red curve), for $\alpha =0.99$ the curve is close to that of the covariance of the ordinary integral of OU (IOU),
given by \eqref{covIOU} (green curve). \par\noindent
In the Figures \ref{surfaceFIOUFISOU} and \ref{2dimFIOUFISOU} of the last section, for further comparisons between the considered processes, we report the three-dimensional plot and two-dimensional color plot, respectively, of the covariance function of FIOU, $c(u,t) = cov({\cal I}^\alpha(Y_{OU})(u), {\cal I}^\alpha(Y_{OU})(t)),$ for various values of $\alpha .$
These illustrate globally the behavior of the
covariance function of FIOU, for any $u$ and $t;$ note, for instance, that for $u=t$ they match the behavior shown in the Figure \ref{varFIOU}.

%
%
%
\bigskip

\subsection{Fractional integrated stationary OU (FISOU) process}
For $\mu >0 $ and $ \sigma \neq 0, $ the {\it stationary OU} (SOU)
process is defined by
\begin{equation} \label{explicitstaOU}
Y_{SOU}(t)= e^{- \mu t} B \left ( \frac { \sigma ^2 } {2 \mu} e^ {2 \mu t } \right );
\end{equation}
therefore, the SOU process is a GM process with
\begin{equation}
r(t)= \frac { \sigma ^2 } {2 \mu} e^ {2 \mu t }, \ m(t)=0,
\end{equation}
\begin{equation}
 h_1(t)= \frac {\sigma ^2 } {2 \mu } e^ { \mu t }, \ h_2(t)= e^ {- \mu t },
 \end{equation}
and covariance
\begin{equation}\label{CovSOU}
c(s,t) = \frac { \sigma ^2 } {2 \mu} e^ {- \mu (t-s)}, \ s \le t .
\end{equation}
One has $var(Y_{SOU}(t))= \frac { \sigma ^2 } {2 \mu}$ and $Y_{SOU}(t)$ admits a
steady-state distribution, which is
\par\noindent
${\cal N} (0, \sigma ^2 / 2 \mu ).$
Moreover,  $r(0)=  \frac {\sigma ^2 } {2 \mu } > 0 ,$ $\eta = B( \frac {\sigma ^2 }
{2 \mu }  ),$  \  $\rho (t)= r(t) - r(0) = \frac {\sigma ^2} {2 \mu} \left (e ^{2
\mu
t } -1 \right ) ,$ $h_1(t)$ and $h_2(t)$ are the same functions  of (non stationary)
OU.

\ep{Referring to the neuronal model \eqref{1}-\eqref{2},  the ${\upeta}(t)$ process solving the equation \eqref{2}, with $\upeta_0$ r.v., is a stationary OU process as 
above specified. A specific example of application can be found in \cite{Sak99}. It is usually used for modeling correlated exogenous inputs.}

\begin{figure}
\centering
\includegraphics[height=0.33 \textheight]{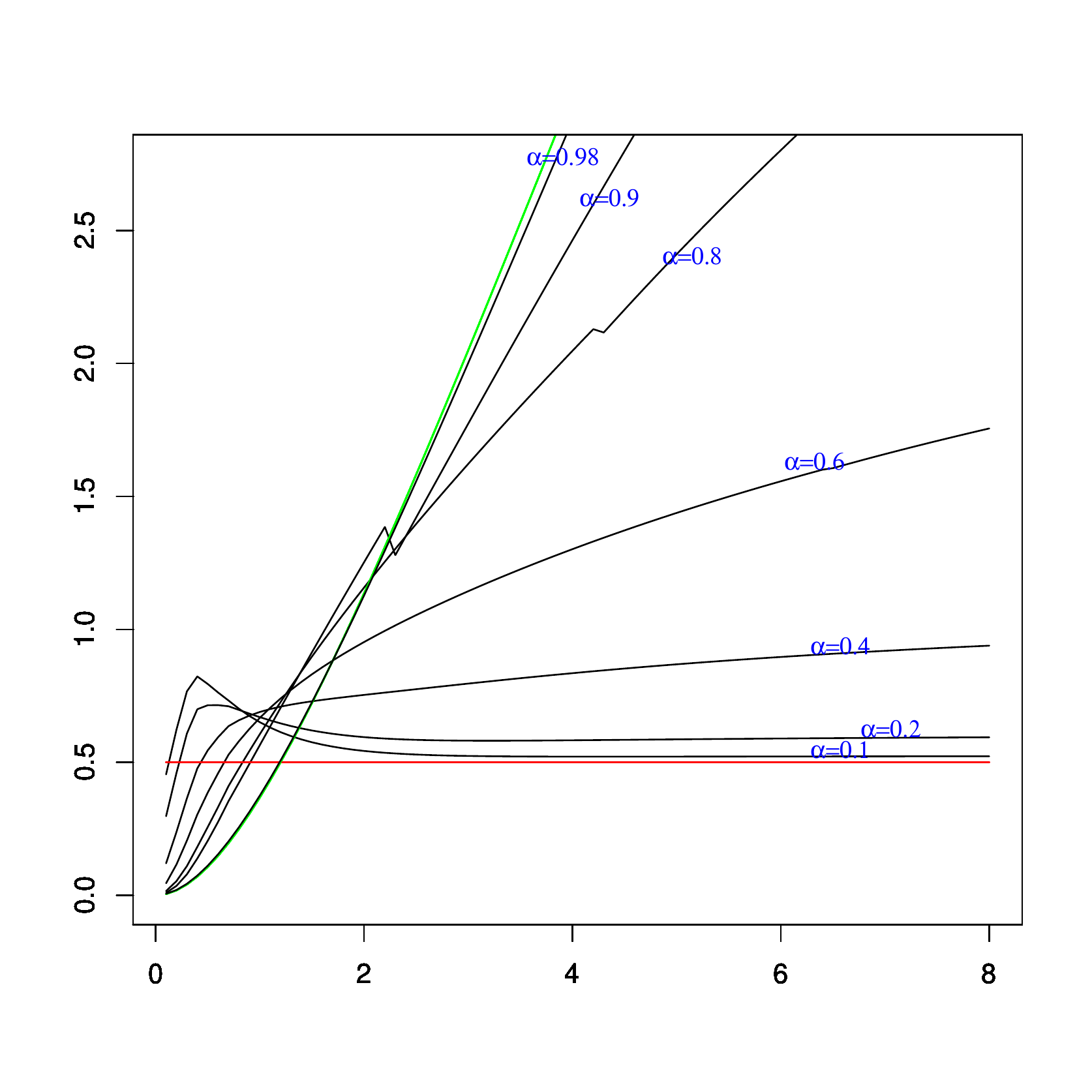}
\caption{Plot of the variance of FISOU, ${\cal I} ^ \alpha (Y_{SOU}) (t)$ with $\sigma = \mu =1,$ as a function of $t \in [0,8],$
for various values of $\alpha \in (0,1).$  The time zone at which an inversion of behavior is observed goes from
$t\simeq  1.7$ to $t \simeq 1.8 \ .$ For $\alpha =0.1$ the curve is close to the graph of the function
$v(t)= \sigma ^2 / 2 \mu$ (i.e. the variance of SOU, red), since FISOU for $\alpha =0$
becomes SOU; for $\alpha \approx 1$ the curve matches the graph of the function $v(t)= \frac { \sigma ^2 } { \mu ^3 } (\mu t + e^ {- \mu t } -1)$
(i.e. the variance of ISOU, green),
given by \eqref{varstaOU}, since  FISOU
becomes the ordinary integral of SOU.
}
\label{varFISOUfig}
\end{figure}

\begin{figure}
\centering
{\includegraphics[width=0.4\textwidth]{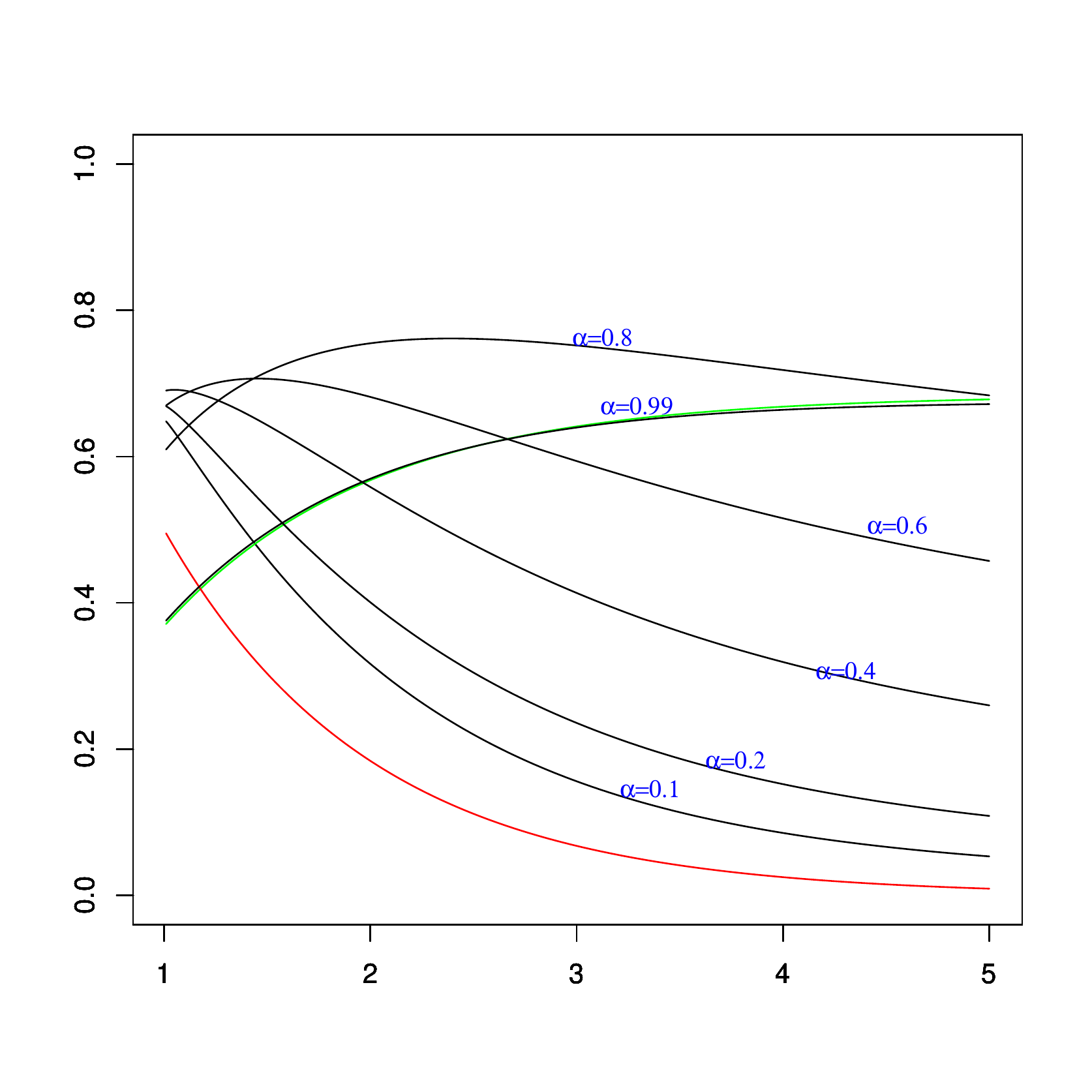}}
\caption{ Plot of the covariance of FISOU, $cov({\cal I}^\alpha(Y_{SOU})(u), {\cal I}^\alpha(Y_{SOU})(t)),$ with $\sigma = \mu =1,$  for  $1=u \le t \le 5,$
and various values of $\alpha \in (0,1).$
For $\alpha = 0.0$, the red curve is the covariance of SOU, given by
\eqref{CovSOU}, for $\alpha =1$ the green curve is the covariance of the ordinary integral of SOU (ISOU),
given by \eqref{CovISOU}.}
\label{covFISOUfig}
\end{figure}

Then, by Theorem \ref{proposition1} the integrated SOU (ISOU) process
(i.e. the ordinary integral over time of the SOU process) is:
\begin{equation} \label{ISOU}
 X_{ISOU}(t)= \int _0 ^ t Y_{SOU}(s) ds = \int _ 0 ^ t e^ {- \mu s } B( \rho (s)) ds
+ \eta \int _0 ^ t e^ {- \mu s } ds
\end{equation}
 \begin{equation}
= \int _ 0 ^ {\rho (t)}  \frac {e^ {- \mu \rho ^ {-1} (u) }} {\rho ' (\rho ^ {-1} (u))} B( u) du + \left ( \frac {1- e^ {- \mu t} } {
\mu}  \right ) \eta .
\end{equation}
Thus, for fixed $t, \  X_{ISOU}(t)$ turns out to be the sum of two (correlated) Gaussian r.v. $W_1$ and $W_2$ with zero mean; the variance of
$W_1$ can be calculated by Theorem \ref{proposition1}, while $var(W_2)= \sigma ^2 \frac {(1- e^ {- \mu t})^2 } { 2 \mu ^3 }.$
The covariance function of ISOU
 $X_{ISOU}(t)$ was calculated in \cite{abupir:physicaA17}; it holds (see also \cite{tay:94}):
\begin{equation}\label{CovISOU}
 cov (X_{ISOU}(t), X_{ISOU}(s)) =  \frac { \sigma ^2 } { 2 \mu ^3 } \left ( 2 \mu
\min (s, t) + e^ {- \mu s } + e^ {- \mu t } - e^ {- \mu |t-s|} -1 \right ) ,
 \end{equation}
and so
\begin{equation} \label{varstaOU}
var (X_{ISOU}(t)) = \frac { \sigma ^2 } { \mu ^3 } \left ( \mu t + e^ {- \mu t} -1
\right ).
\end{equation}
As far as the fractional integral ${\cal I} ^ \alpha (Y_{SOU}) (t)$  (FISOU) is concerned, by using \eqref{ISOU} one has:
\begin{equation} \label{defFISOU}
{\cal I}^ \alpha  (Y_{SOU}) (t)= {\cal I}^ \alpha (g _1 ) (t) +  \eta \ {\cal I}^ \alpha  (g_2 ) (t),
\end{equation}
with
\begin{equation}
g_1 (t)= e^ {- \mu t} B (\rho (t) ) , \ g_2 (t) = e^ { - \mu t },
\end{equation}
where the distributions of the two fractional integrals can be obtained  by using Proposition \ref{proposizionedue}.
By using \eqref{defFISOU}, we obtain the covariance of FISOU, for $ u \le t :$
$$
cov ({\cal I}^ \alpha  (Y_{SOU}) (u),{\cal I}^ \alpha  (Y_{SOU}) (t)) =
E  \left [ {\cal I}^ \alpha  (e^{- \mu u } B( \rho (u))) \ {\cal I}^ \alpha  (e^ { - \mu t } B( \rho (t)))  \right ] +
$$
$$
E \left [ {\cal I}^ \alpha  (e^{- \mu u } B( \rho (u)))  \ \eta {\cal I} ^ \alpha (e^ { - \mu t } ) \right ]+
E \left [\eta {\cal I} ^ \alpha (e^ { - \mu u } ) \  {\cal I}^ \alpha  (e^{- \mu t } B( \rho (t)) \right ] +
$$
\begin{equation} \label{covFISOU}
E  \left [ \eta {\cal I} ^ \alpha (e^ { - \mu u } ) \  \eta {\cal I} ^ \alpha (e^ { - \mu t } ) \right ] := \widetilde J_1+  \widetilde J_2 + \widetilde J_3 + \widetilde J_4 ,
\end{equation}
where for simplicity we let drop the dependence of $\widetilde J _i $ on $ \alpha .$
By using the definition of fractional integral, and proceeding as in the case of FIBM and FIOU, one gets:
$$ \widetilde J_1 = \frac {\widetilde J _ {11} + \widetilde J_ {12} + \widetilde J _{13} } {\Gamma ^2 ( \alpha ) } ,$$
where:
$$ \widetilde J _ {11} = \frac {\sigma ^2 } { 2 \mu } \int _ 0 ^u ds (u-s)^ { \alpha -1} e ^ { - \mu s} \int _ 0 ^s dv (t-v) ^ {\alpha -1}
e^ { - \mu v } \left (e ^{2 \mu v } -1 \right ) = \widetilde I_1, $$
$$ \widetilde J _ {12} = \frac {\sigma ^2 } { 2 \mu } \int _ 0 ^u ds (u-s)^ { \alpha -1} e ^ { - \mu s}  \left (e ^{2 \mu s } -1 \right )\int _ s ^u dv (t-v) ^ {\alpha -1} e ^ { - \mu v } = \widetilde I _2, $$
$$ \widetilde J _ {13} = \frac {\sigma ^2 } { 2 \mu }  \left ( \int _ 0 ^u ds (u-s)^ { \alpha -1} e ^ { - \mu s}  \left (e ^{2 \mu s } -1 \right )  \right )
\left ( \int _ u ^t dv (t-v) ^ {\alpha -1} e ^ { - \mu v } \right ) = \widetilde I_3 ,$$
being $\widetilde I _i$ the integrals concerning FIOU, given by \eqref{integraliItilde}.
Also:
$$ \widetilde J _ {2} =   \frac { \sigma ^2 } {2 \mu }  \left ( \int _ 0 ^u ds (u-s) ^ {\alpha -1 } e ^ { - \mu s }
\cdot \min  \large \{ 1, e ^{2 \mu s } -1 \large  \} \right )
\left ( \int _ 0 ^t dv (t-v) ^ {\alpha -1} e ^ { - \mu v } \right ),$$
$$ \widetilde J _ {3} =  \widetilde J _ {2} , $$
$$
\widetilde J _ {4}=  \frac { \sigma ^2 } {2 \mu }  \left ( \int _ 0 ^u ds (u-s) ^ {\alpha -1 } e ^ { - \mu s }
 \right )
\left ( \int _ 0 ^t dv (t-v) ^ {\alpha -1} e ^ { - \mu v } \right ).$$
Finally, substituting in \eqref{covFISOU}, we obtain for $u \le t$ the covariance of FISOU:
$$
cov ({\cal I}^ \alpha  (Y_{SOU}) (u),{\cal I}^ \alpha  (Y_{SOU}) (t))   =
\frac 1 { \Gamma ^2 ( \alpha )} \left ( \widetilde J _{11} + \widetilde J _{12} + \widetilde J _{13} + 2 \widetilde J _{2} + \widetilde J_4 \right )
$$
\begin{equation} \label{covFISOUbis}
= \frac 1 { \Gamma ^2 ( \alpha )} \left ( \widetilde I_1 + \widetilde I_2 + \widetilde I_3 + 2 \widetilde J _{2}  + \widetilde J _4 \right ) =
cov ({\cal I}^ \alpha  (Y_{OU}) (u),{\cal I}^ \alpha  (Y_{OU}) (t)) + \frac {2 \widetilde J _{2}  + \widetilde J _4} { \Gamma ^2 ( \alpha) } .
\end{equation}
Thus, the covariance function of FISOU turns out to be the sum of the covariance of FIOU and the quantity $\frac {2 \widetilde J _{2}  + \widetilde J _4} { \Gamma ^2 ( \alpha) } .$
\par\noindent
For $u=t, $ one has $\widetilde J _{13} =0$ and $\widetilde J _{11} = \widetilde J _{12},$ so:
\begin{equation} \label{varFISOU}
var ({\cal I} ^ \alpha ( Y_{SOU} (t) ) = \frac 1 { \Gamma ^ 2 ( \alpha ) } \left ( 2 \widetilde J _{11} + 2 \widetilde J _{2}  + \widetilde J _4 \right )
= \frac 1 { \Gamma ^ 2 ( \alpha ) } \left ( 2 \widetilde I_1 + 2 \widetilde J _{2}  + \widetilde J _4 \right ),
\end{equation}
where all the integrals have to be calculated for $u=t.$
 \par\noindent

In the Figures \ref{varFISOUfig}-\ref{covFISOUfig},
we report the shape of the variance and covariance function of FISOU,
for various values of $\alpha \in (0,1).$  In the Figure \ref{varFISOUfig}, it is observed a time zone (from $t \simeq 1.7$ to $t \simeq 1.8$), at which there is an inversion of behavior, as in Figure \ref{varFIOU} concerning the variance of FIOU. Here, however, unlike Figure \ref{varFIOU}, it is visible a more confused behavior in the time zone $(1.7, 1.8).$ This is probably due to an higher degree of stochasticity of FISOU with respect to FIOU (in fact, the starting point of SOU is random).
Note that we have not be able to calculate the variance and covariance 
of FISOU for $\alpha$ less than $0.1,$ owing to numerical problems which 
arise in computing  the integrals $ \tilde J,$ for $\alpha$ near zero.
For $\alpha =0.1$ the curve is close to the graph of the function
$v(t)= \sigma ^2 / 2 \mu$ (i.e. the variance of SOU), since FISOU for $\alpha =0$
becomes SOU; for $\alpha \approx 1$ the curve matches the graph of the function $v(t)= \frac { \sigma ^2 } { \mu ^3 } (\mu t + e^ {- \mu t } -1)$
(i.e. the variance of ISOU),
given by \eqref{varstaOU}, since  FISOU
becomes the ordinary integral of SOU.

\section{{Some graphical comparisons  and concluding remarks}}

In order to complete our analysis and compare the obtained results for the considered processes, we provide the numerical evaluations of the covariance functions in some three and two dimensional plots.

\begin{figure}[h]
\centering
{\includegraphics[width=0.32\textwidth]{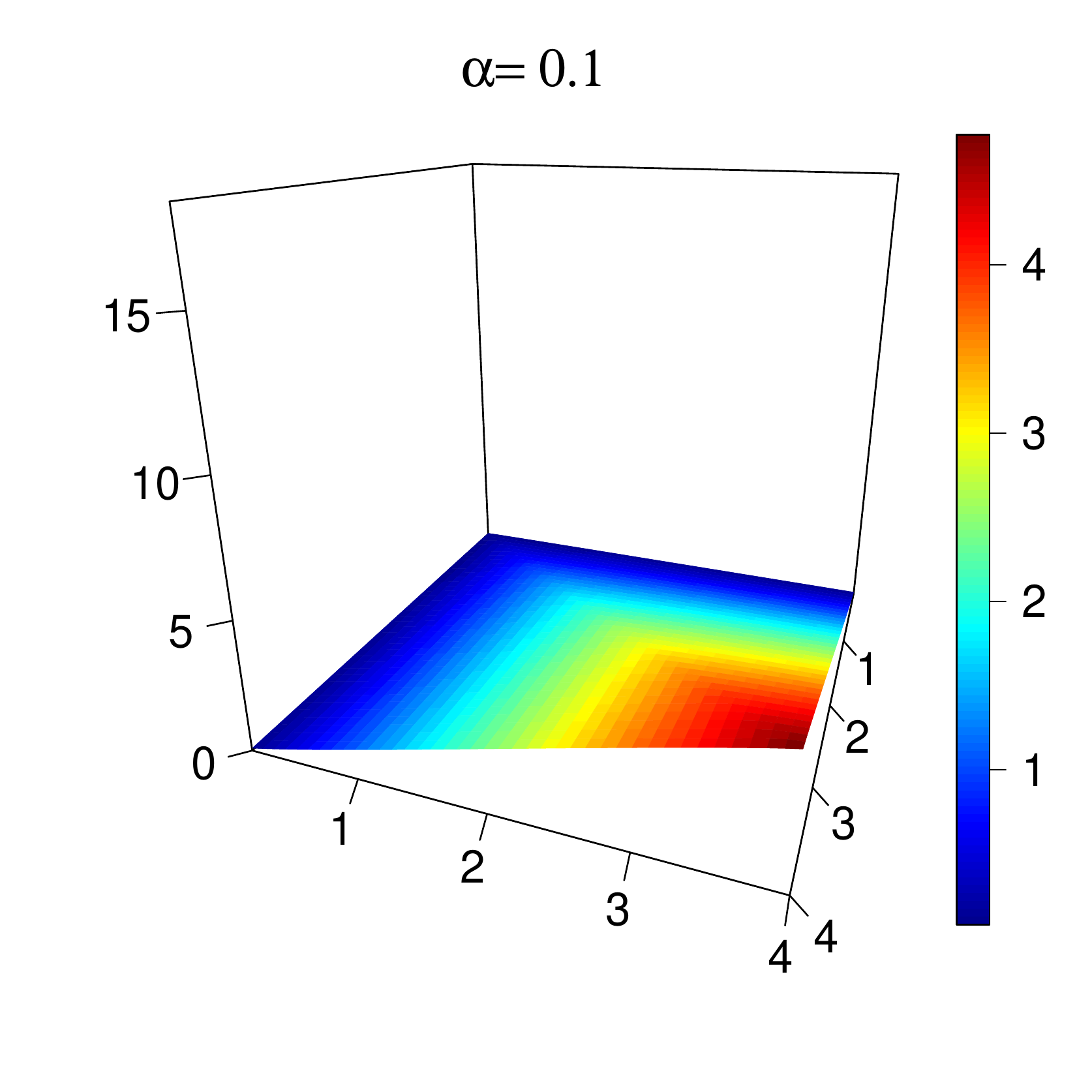}}
{\includegraphics[width=0.32\textwidth]{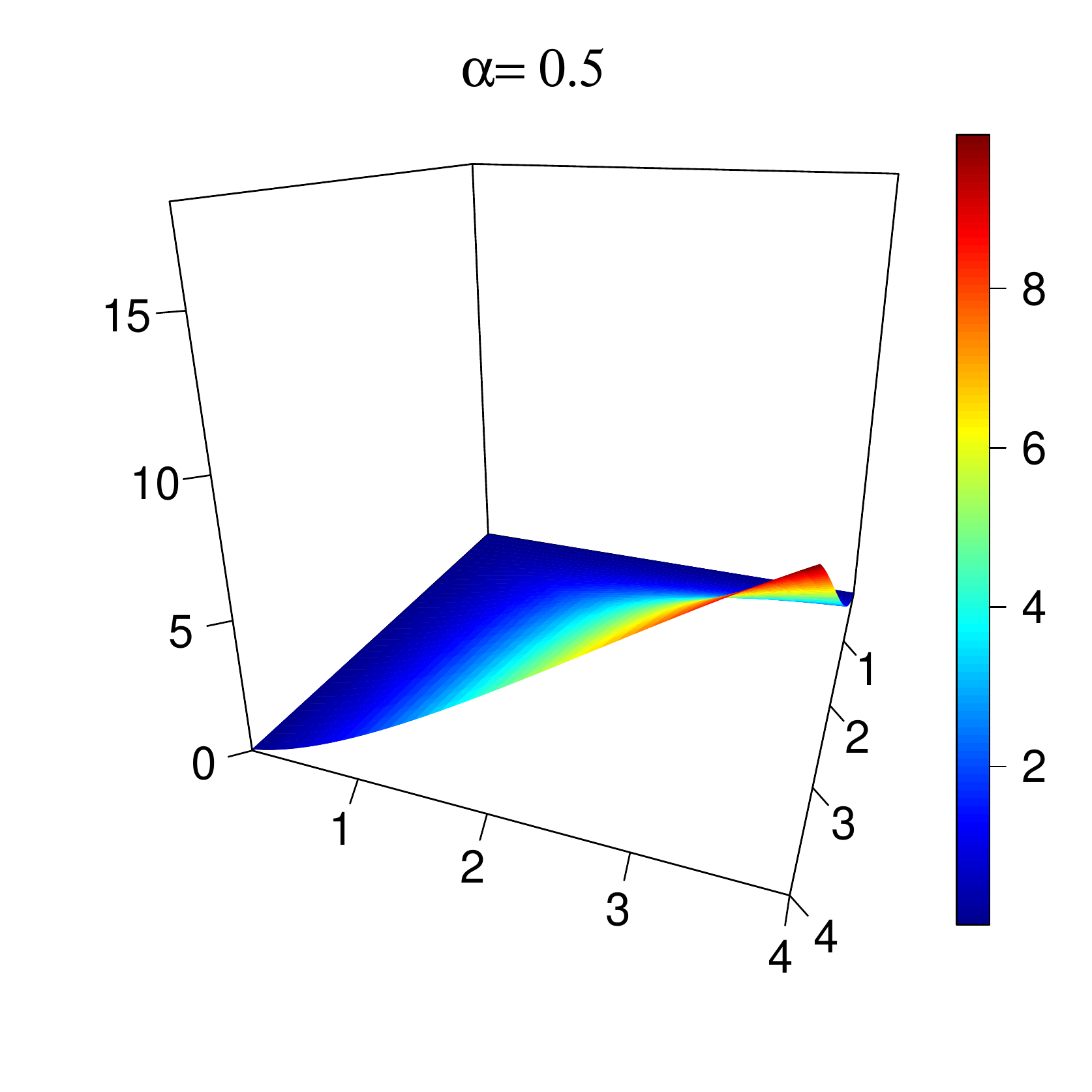}}
{\includegraphics[width=0.32\textwidth]{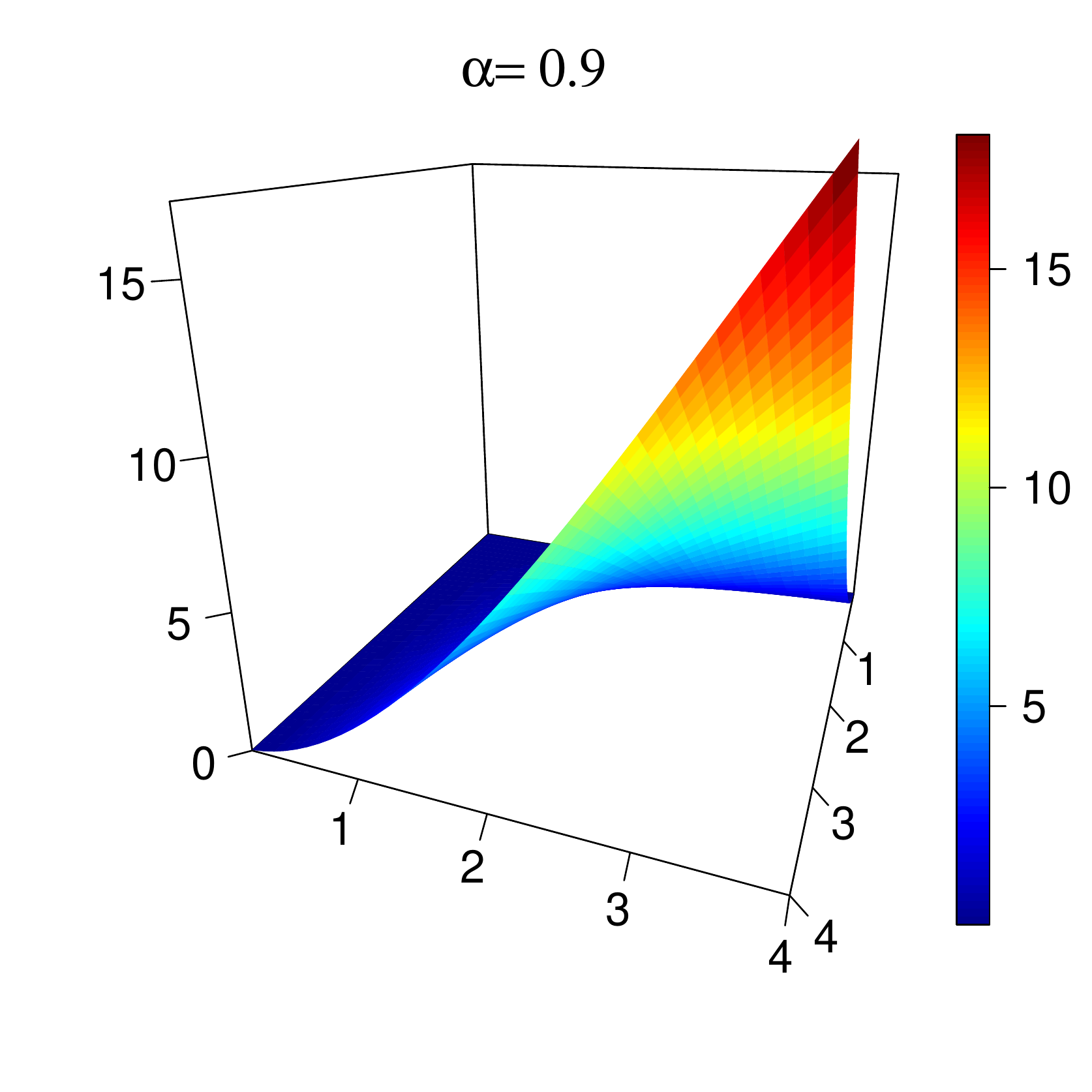}}\\
{\includegraphics[width=0.32\textwidth]{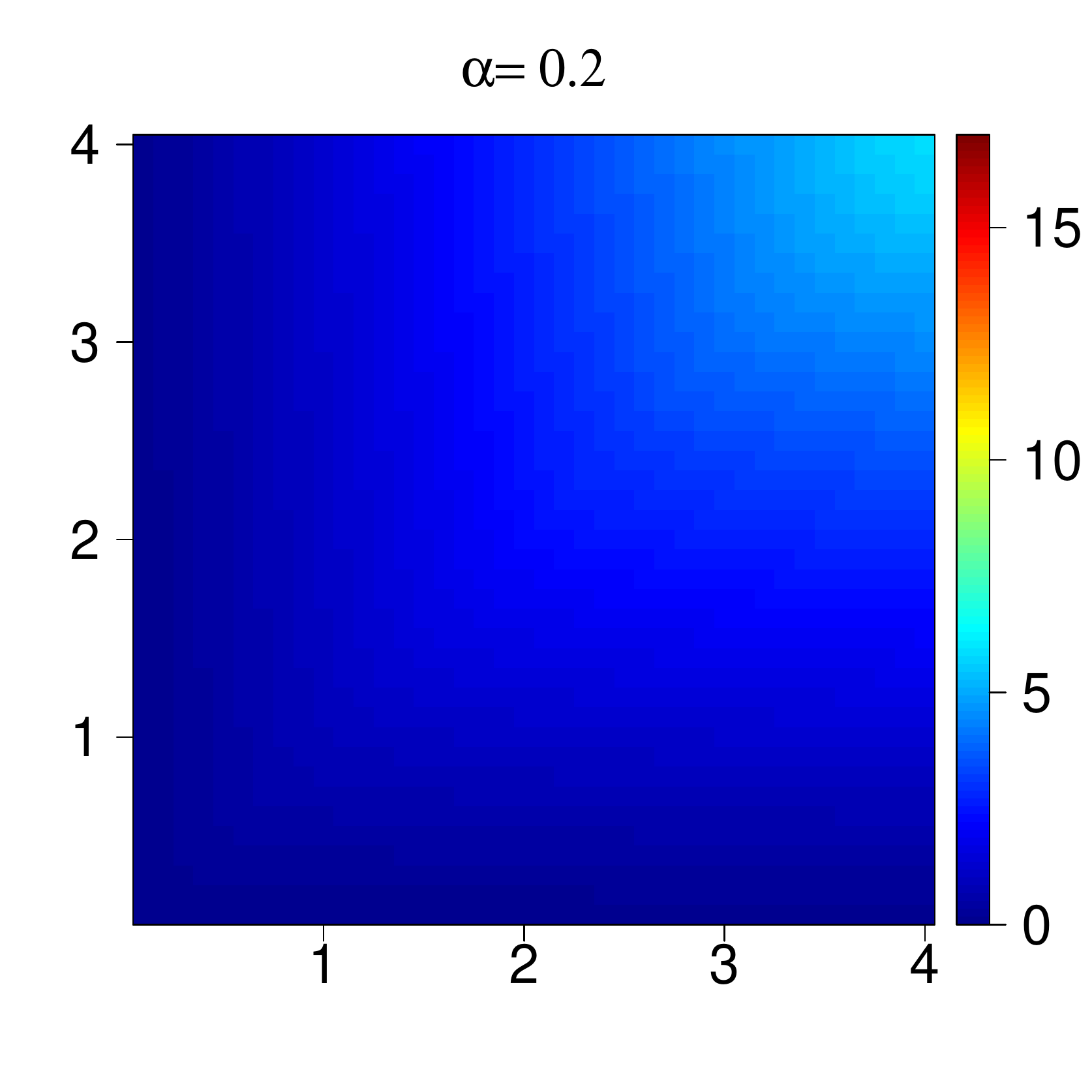}}
{\includegraphics[width=0.32\textwidth]{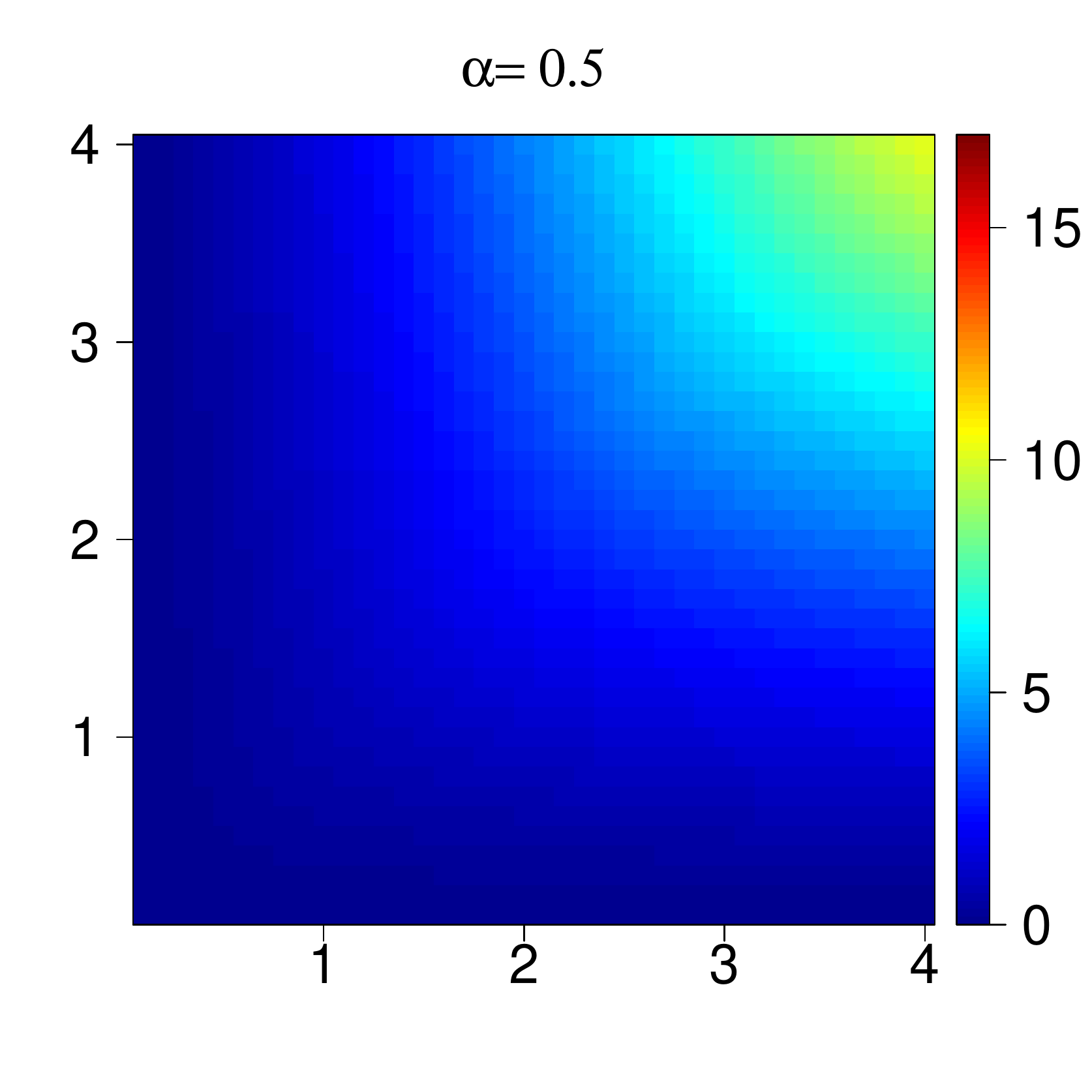}}
{\includegraphics[width=0.32\textwidth]{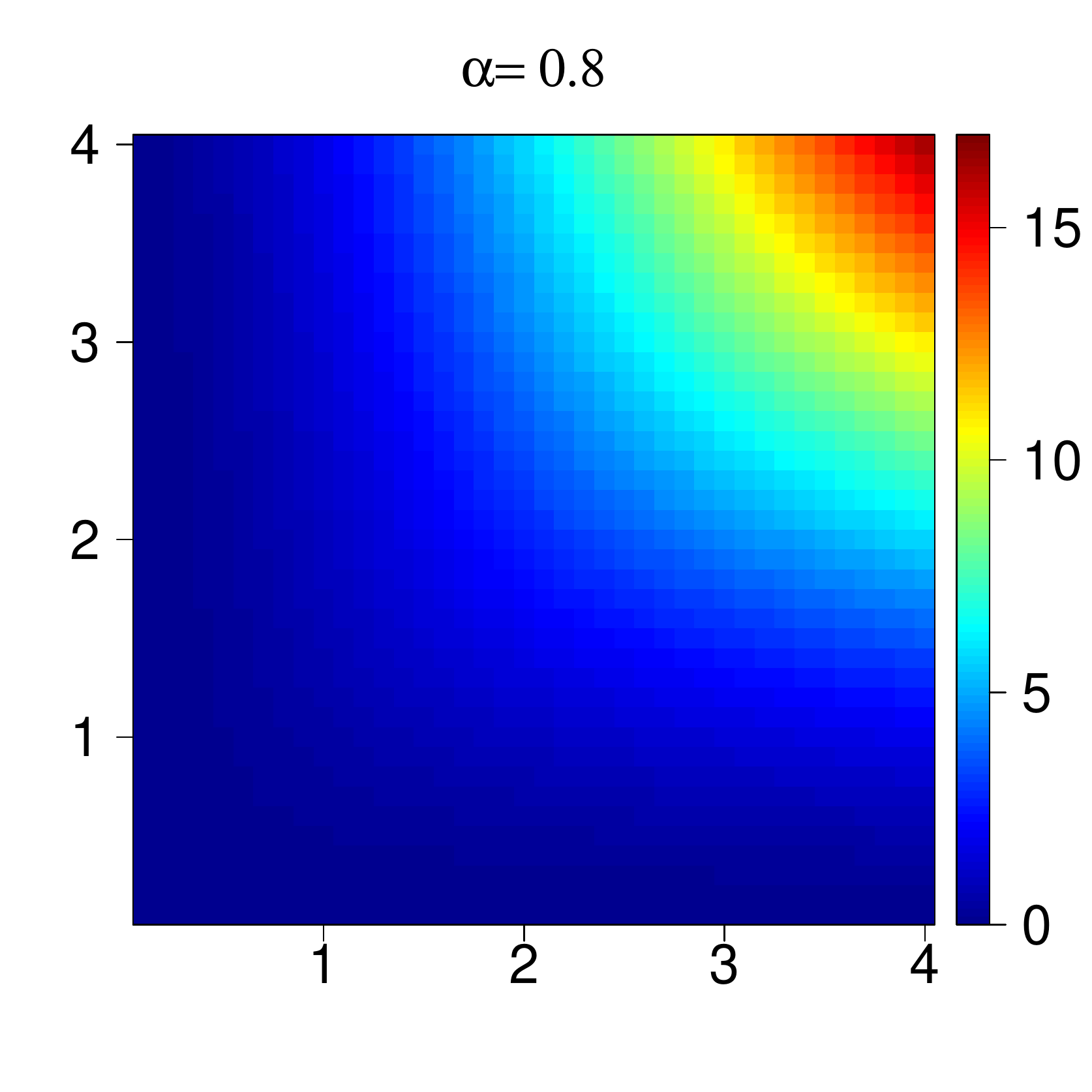}}
\caption{On top: three-dimensional plot of the covariance function of FIBM $cov({\cal I }^ \alpha  (B) (u), {\cal I }^ \alpha  (B)(t))$,  for
$ \sigma = \mu =1,$ and  various values of $\alpha .$
On bottom: two-dimensional color plot of the covariance function of FIBM $cov({\cal I }^ \alpha  (B) (u), {\cal I }^ \alpha  (B)(t))$ on the same scale.}
\label{2dimcolorcovFIBM}
\end{figure}

\begin{figure}[h]
\centering
{\includegraphics[width=0.32\textwidth]{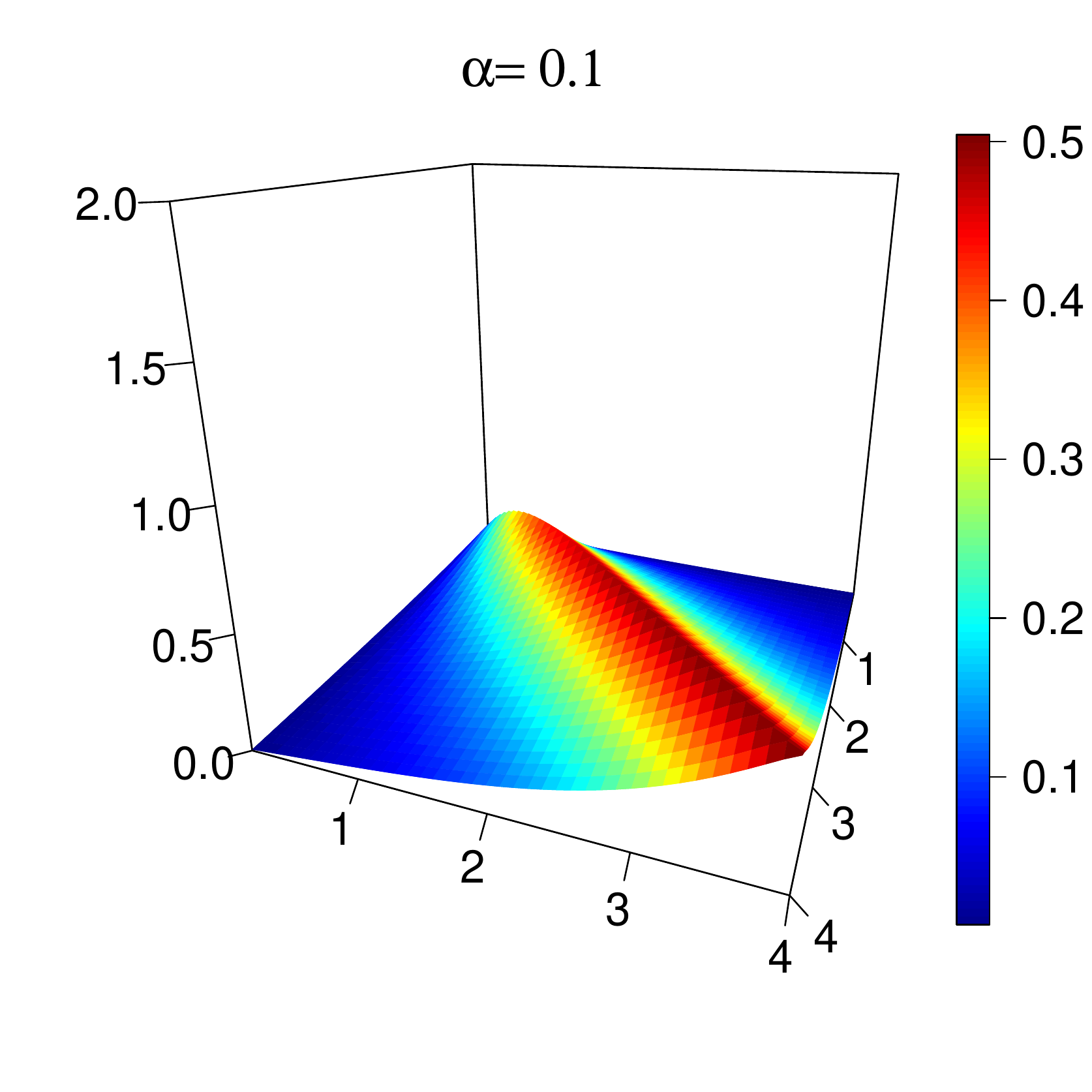}}
{\includegraphics[width=0.32\textwidth]{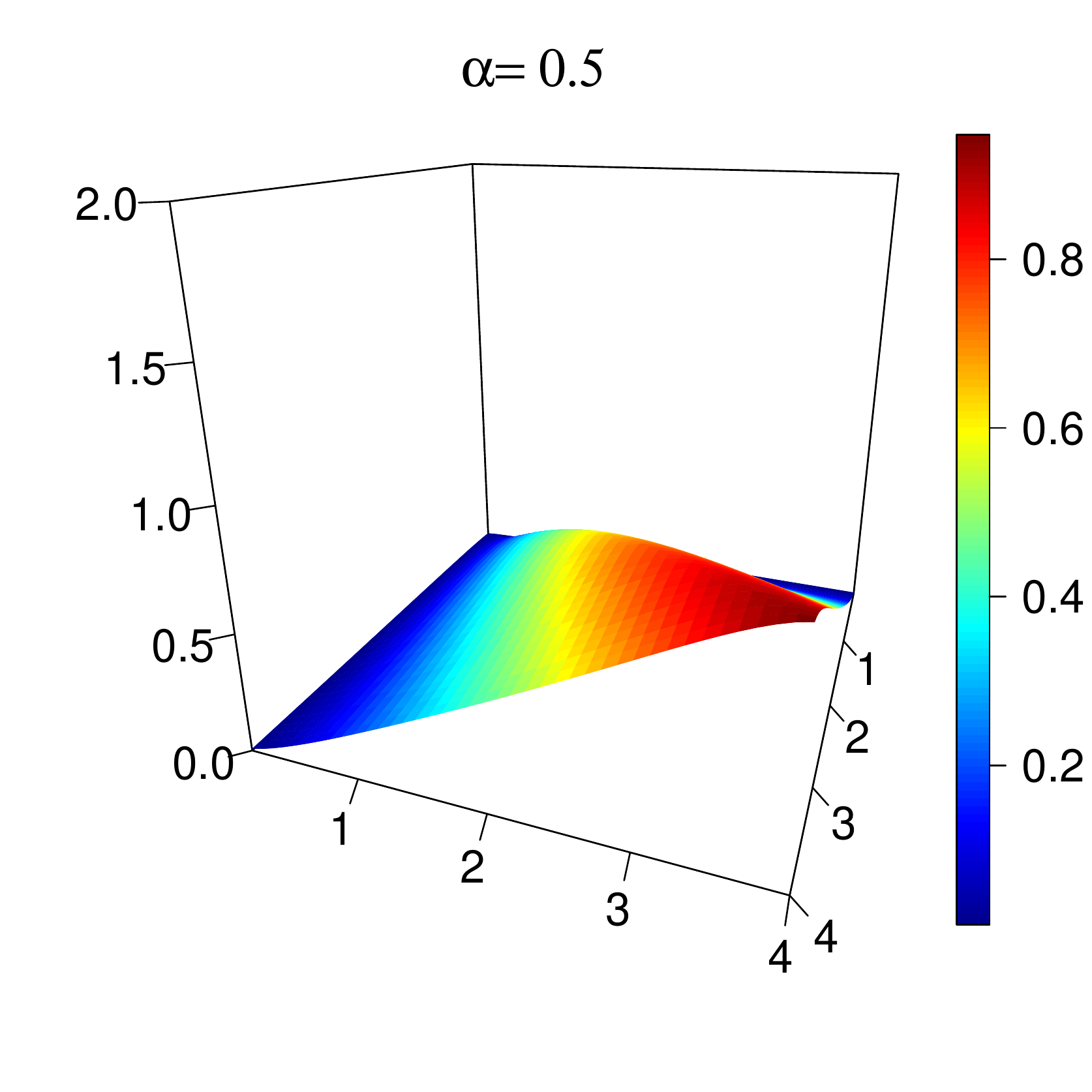}}
{\includegraphics[width=0.32\textwidth]{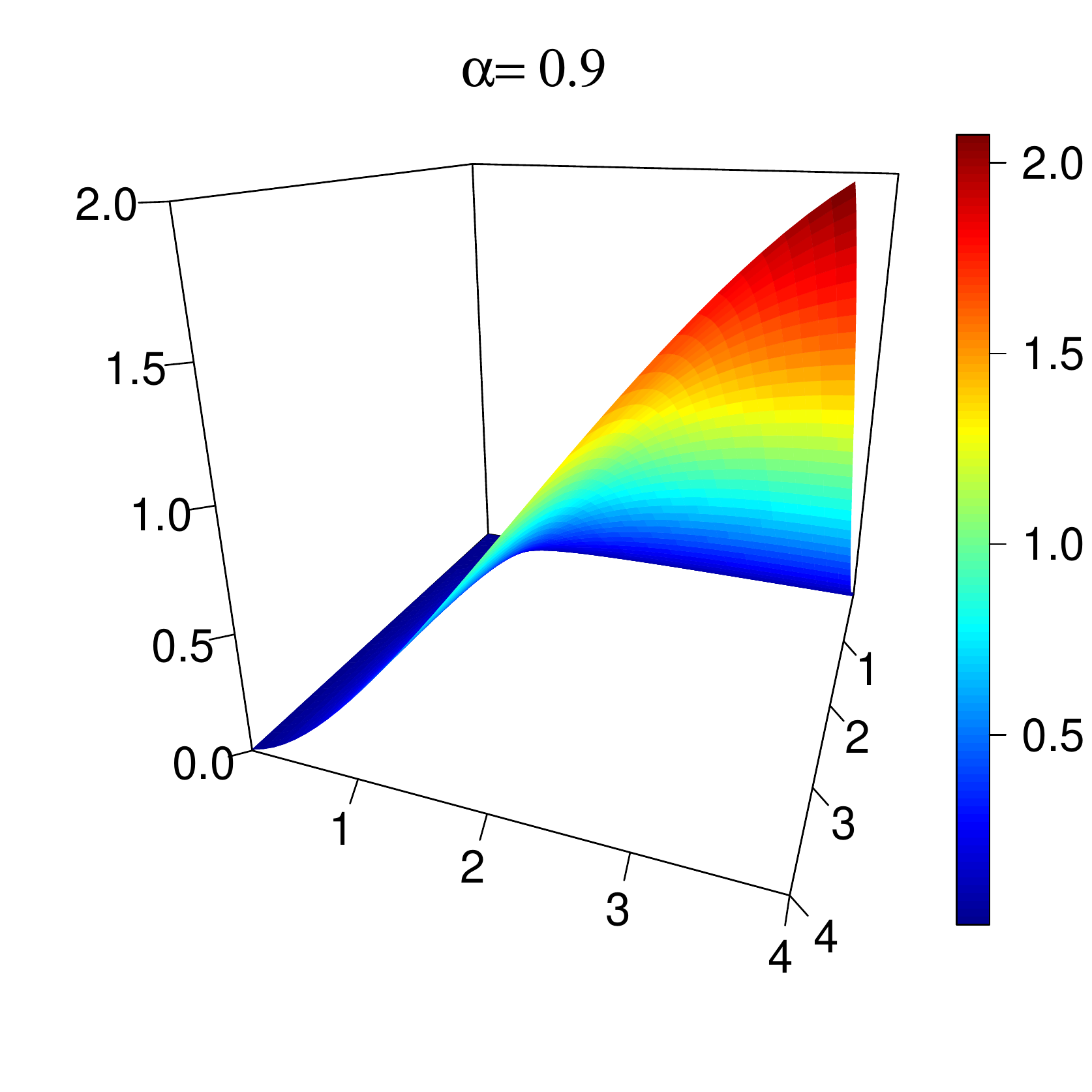}}\\
{\includegraphics[width=0.32\textwidth]{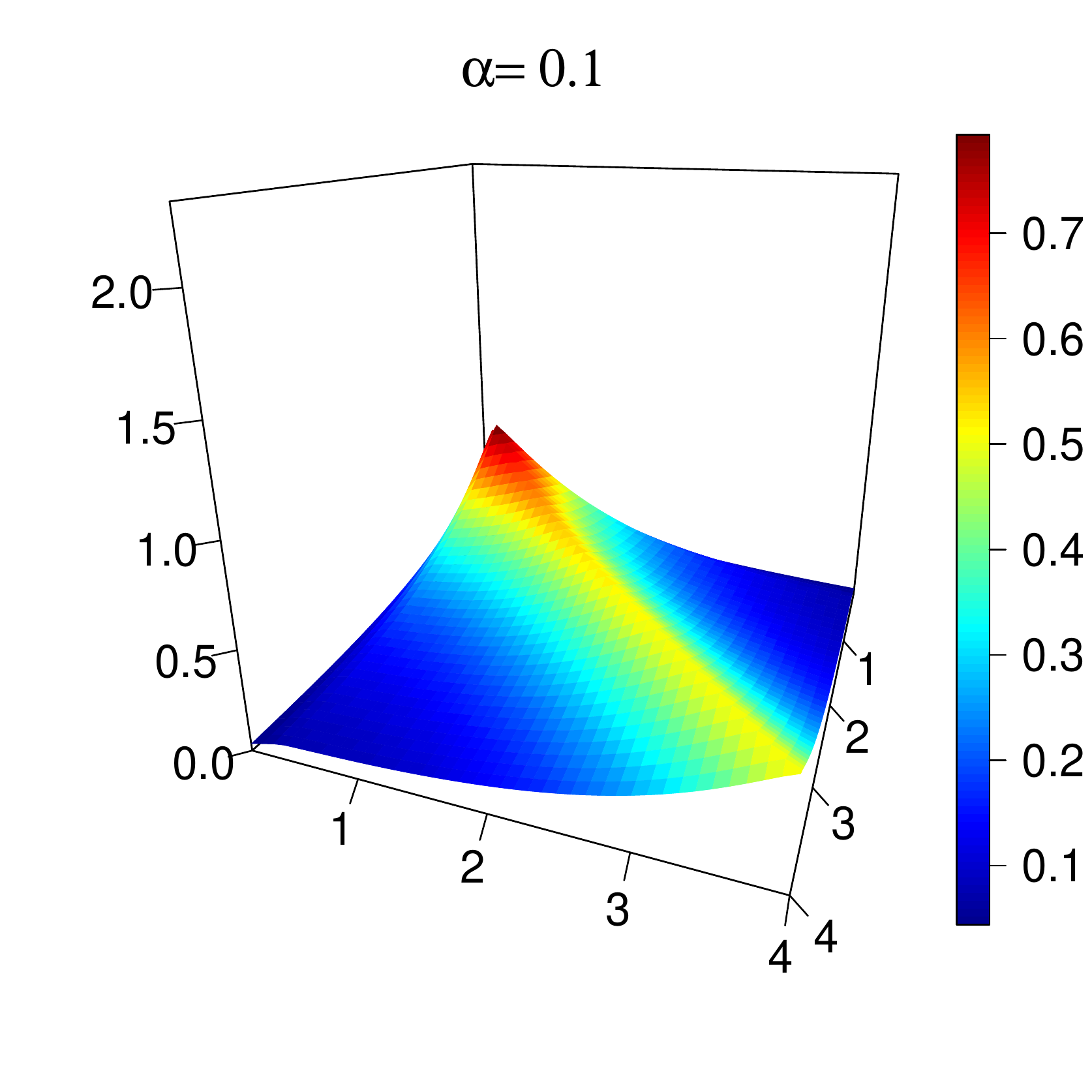}}
{\includegraphics[width=0.32\textwidth]{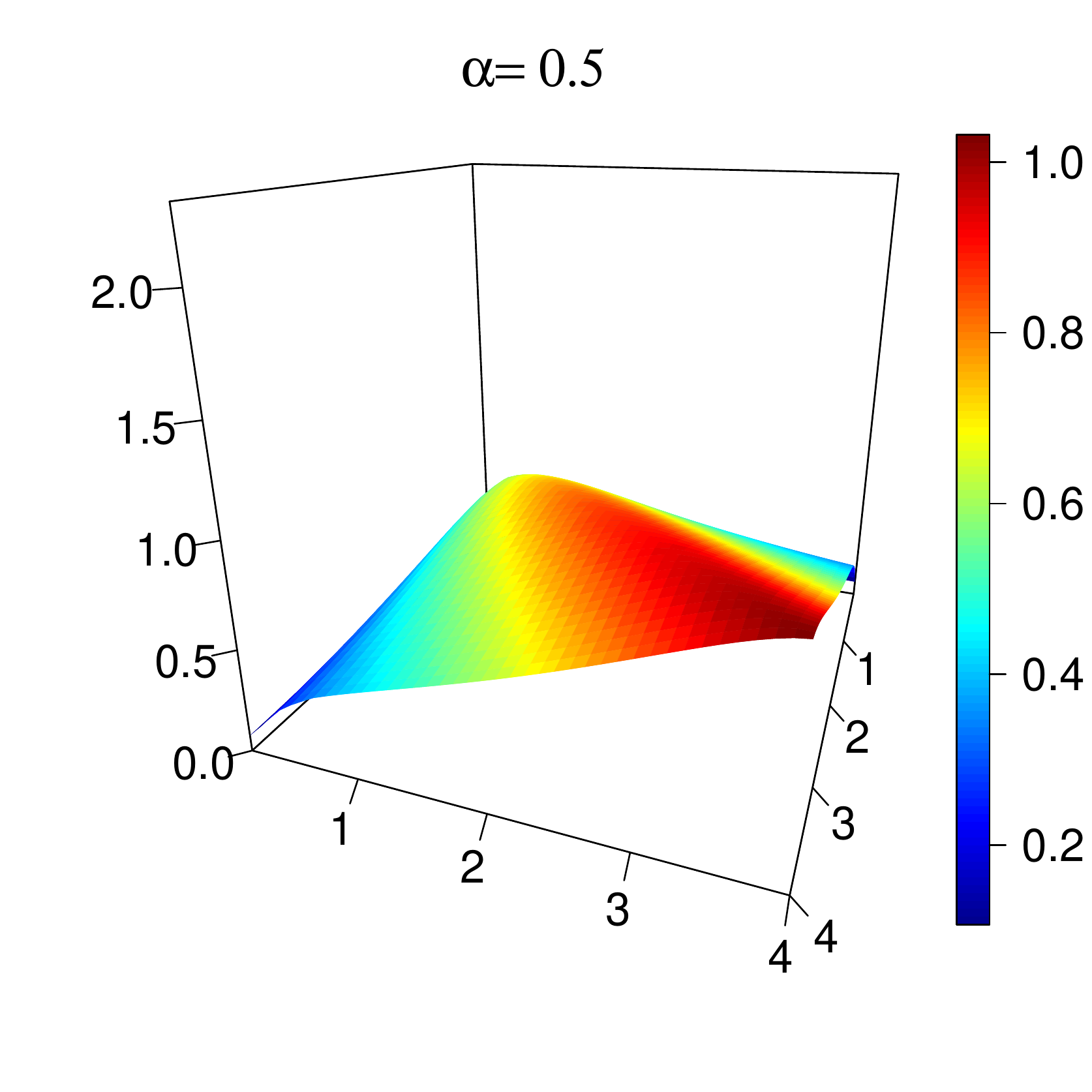}}
{\includegraphics[width=0.32\textwidth]{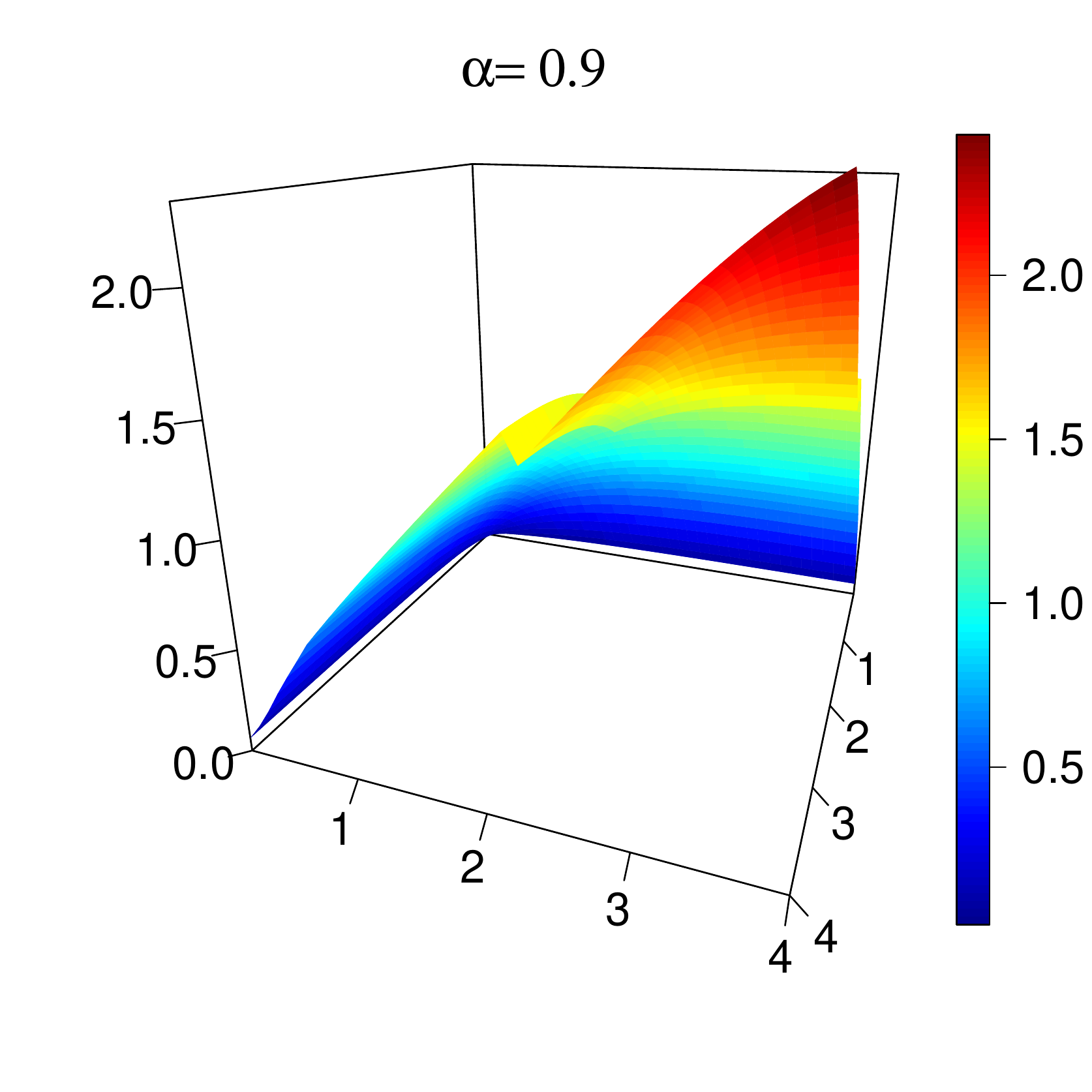}}
\caption{On top: three-dimensional plot of the covariance function of FIOU $c(u,t) = cov({\cal I}^\alpha(Y_{OU})(u), {\cal I}^\alpha(Y_{OU})(t))$  for
$ \sigma = \mu =1,$ and  various values of $\alpha.$ On bottom: the same for the covariance function of FISOU $c(u,t) = cov({\cal I}^\alpha(Y_{SOU})(u), {\cal I}^\alpha(Y_{SOU})(t))$.}
\label{surfaceFIOUFISOU}
\end{figure}

\begin{figure}[h]
\centering
{\includegraphics[width=0.32\textwidth]{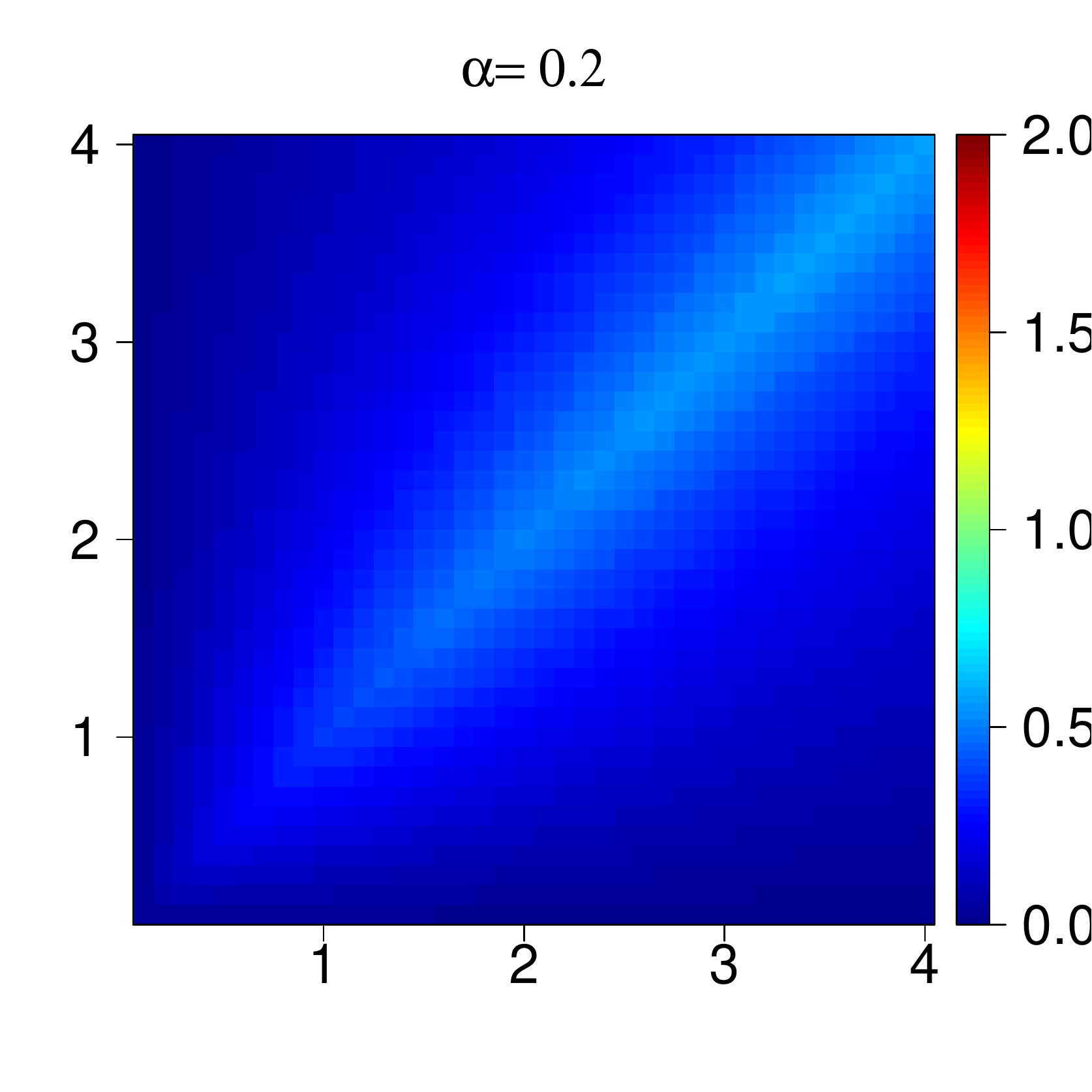}}
{\includegraphics[width=0.32\textwidth]{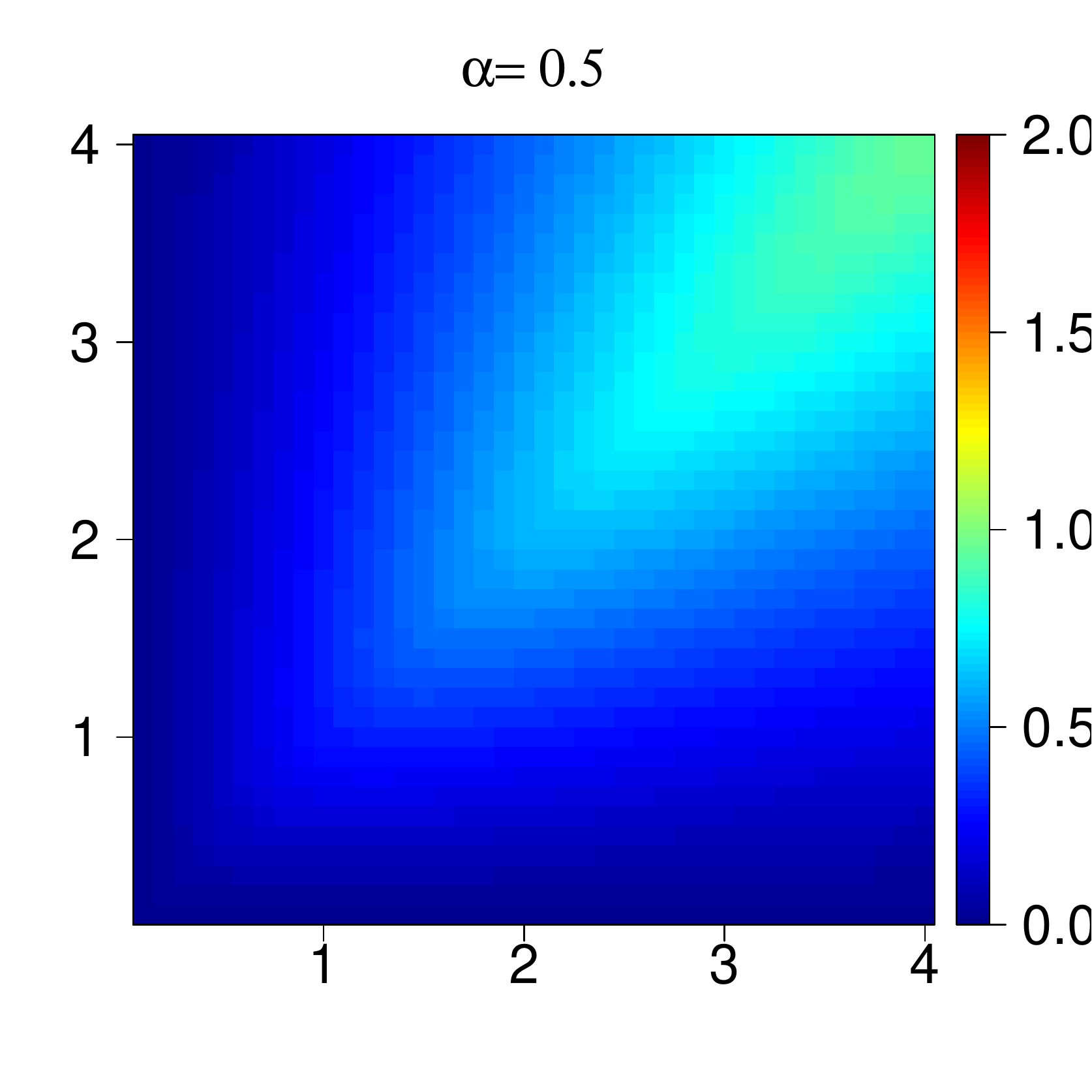}}
{\includegraphics[width=0.32\textwidth]{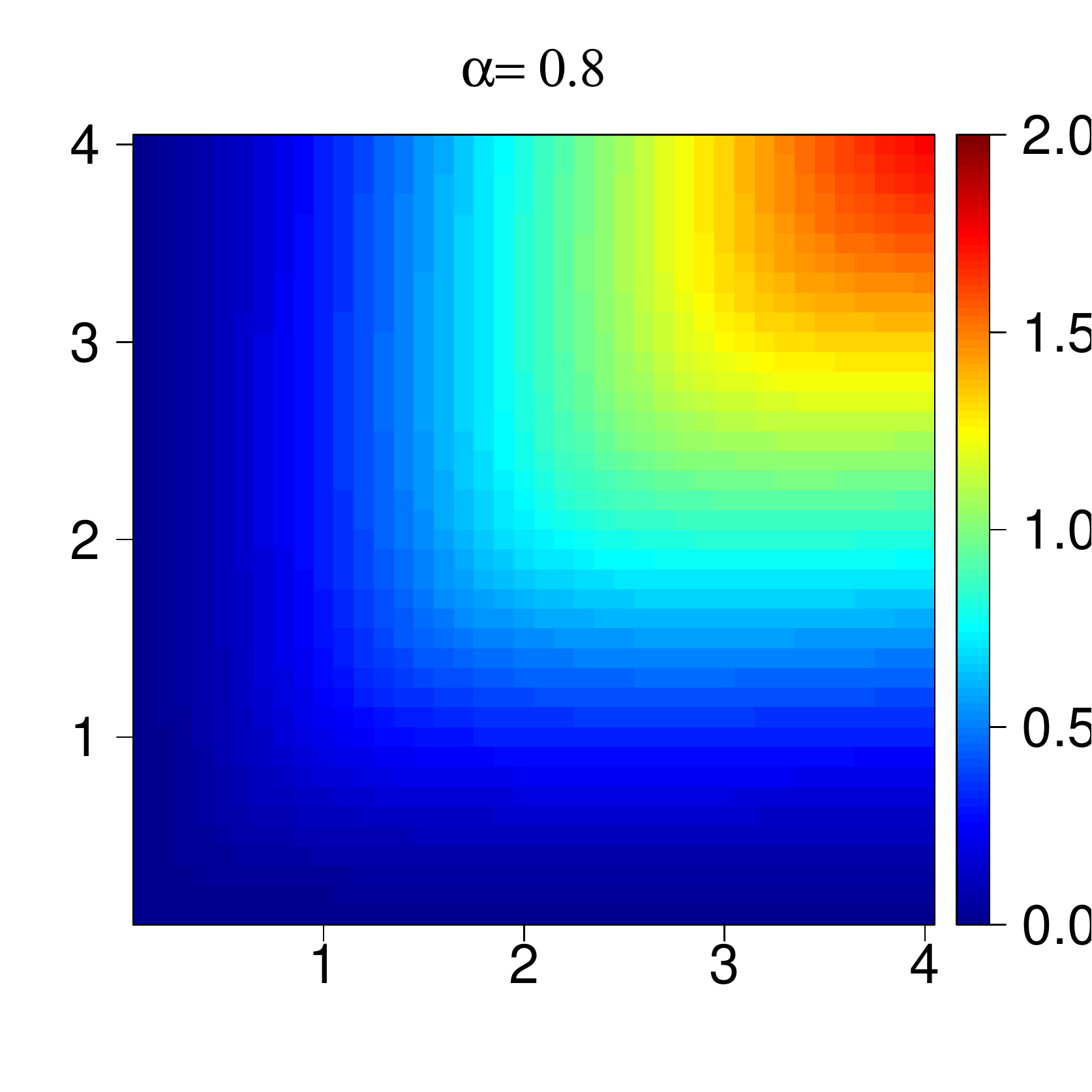}}\\
{\includegraphics[width=0.32\textwidth]{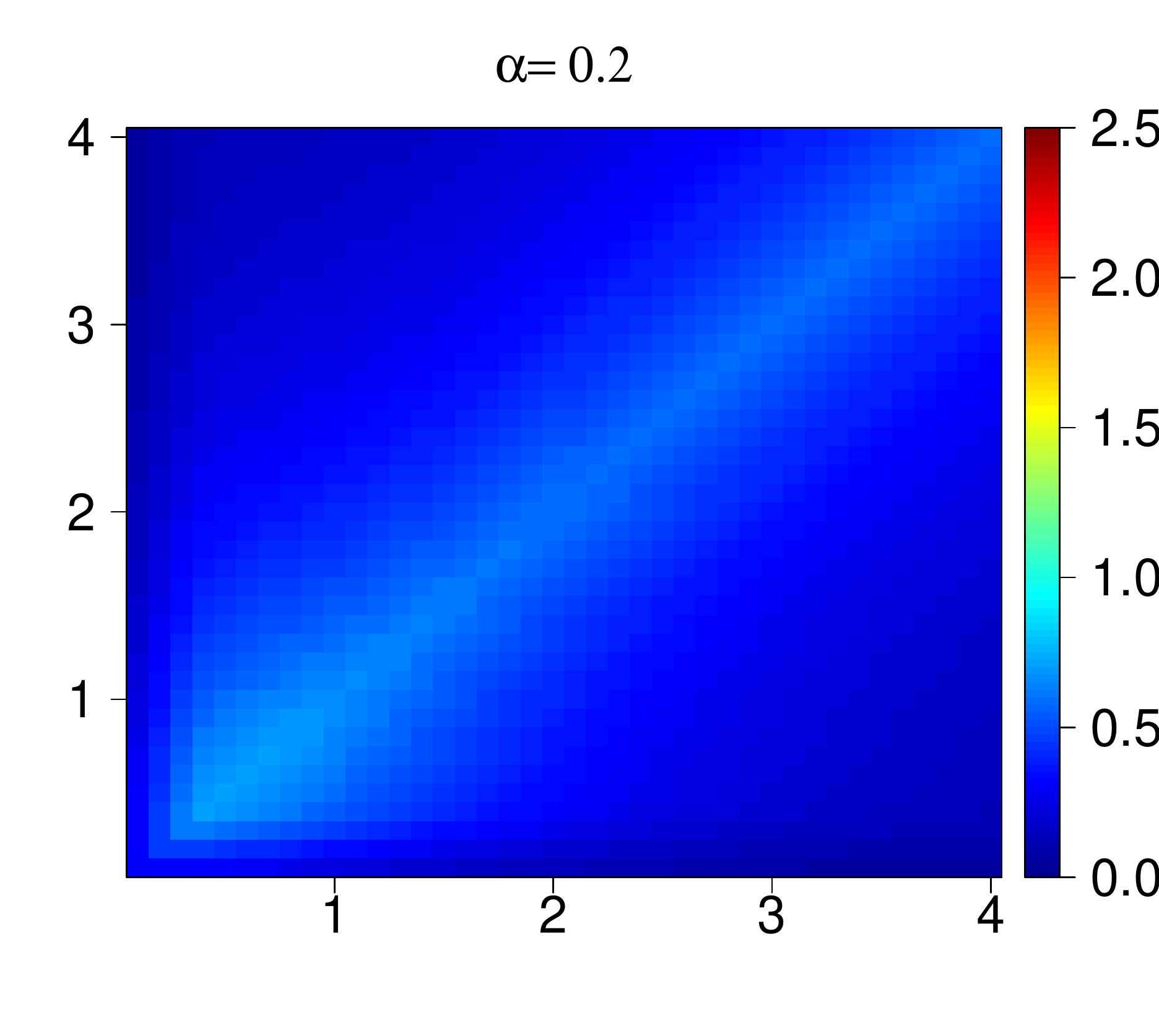}}
{\includegraphics[width=0.32\textwidth]{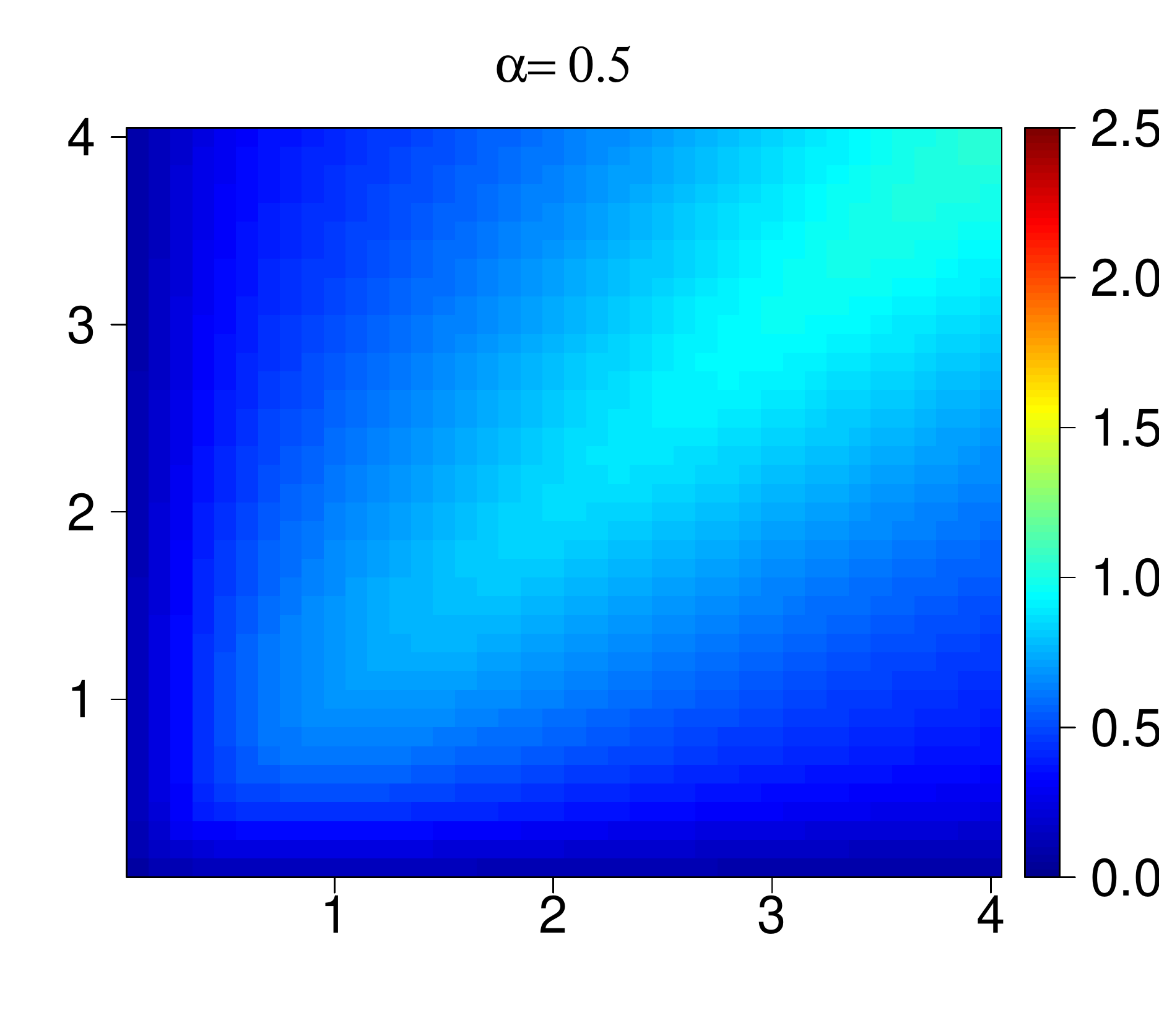}}
{\includegraphics[width=0.32\textwidth]{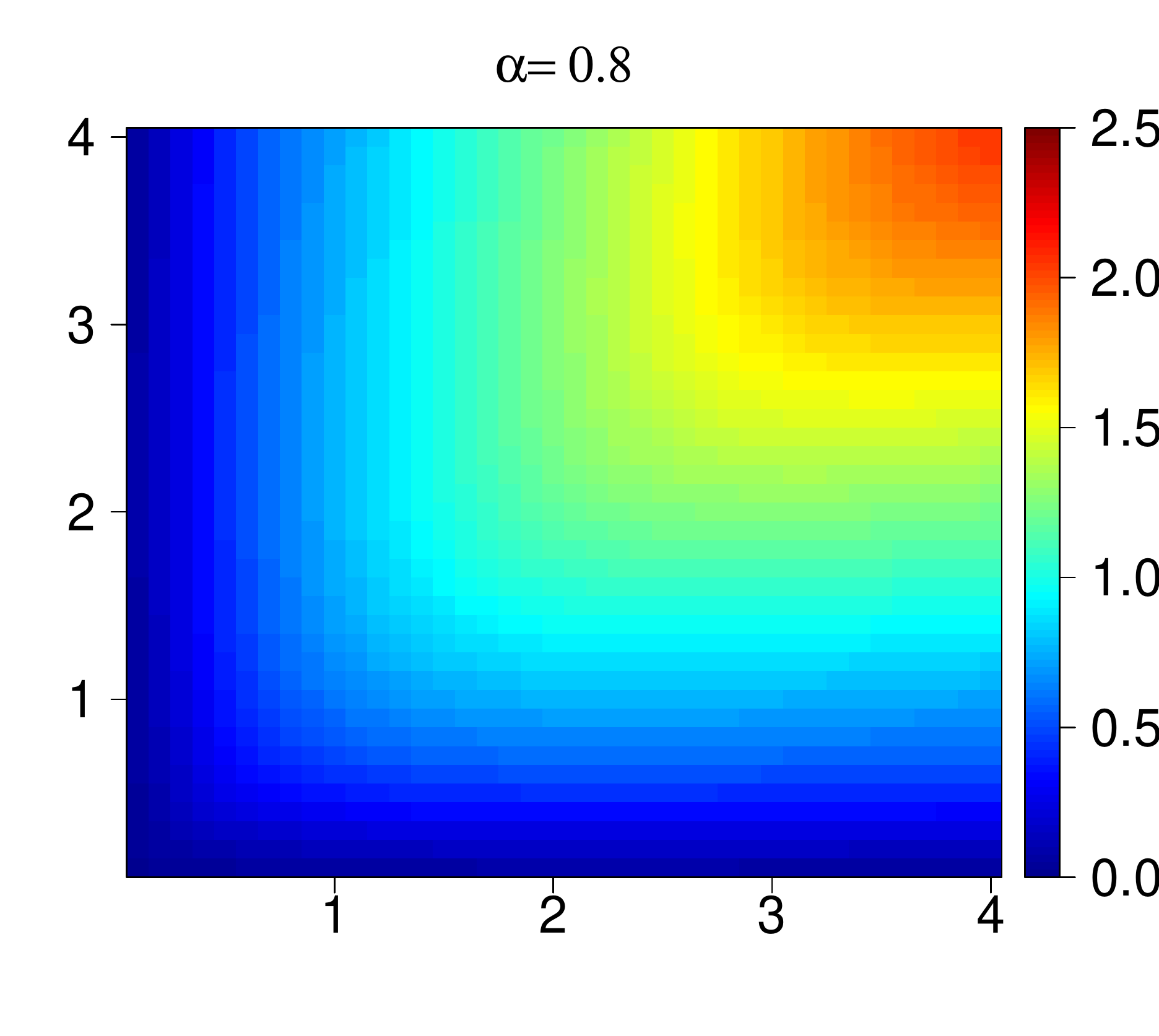}}
\caption{On top: two-dimensional color plot of the covariance function of FIOU $c(u,t) = cov({\cal I}^\alpha(Y_{OU})(u), {\cal I}^\alpha(Y_{OU})(t)),$ for
$ \sigma = \mu =1,$ (on the same scale) and
 various values of $\alpha.$ On bottom: the same for the covariance function of FISOU $c(u,t) = cov({\cal I}^\alpha(Y_{SOU})(u), {\cal I}^\alpha(Y_{SOU})(t)).$}
\label{2dimFIOUFISOU}
\end{figure}

In the Figure\ref{2dimcolorcovFIBM}  we report the two-dimensional color plot of the 
covariance function of FIBM, $c(u,t) = cov({\cal I}^\alpha(B)(u), {\cal 
I}^\alpha(B)(t)),$ given by \eqref{covarianceofI(t)}  and the graph of 
the surface $z = c(u, t),$ for various values of $\alpha ;$  in the 
Figure \ref{surfaceFIOUFISOU} we report
the graphs of the surface $z = c(u, t),$ for various values of $\alpha 
,$ in the cases of FIOU and FISOU, respectively. In the Figure \ref{2dimFIOUFISOU}, we 
report
the two-dimensional color plots of the covariance function of FIOU and 
FISOU, respectively, for various values of $\alpha .$
These figures illustrate globally the behavior of the
covariance function of the three considered processes, for any $u$ and 
$t;$ note, for instance, that for $u=t$ they match the behavior shown in 
the figures concerning the variance of the processes (compare with Figures \ref{varFIOU},\ref{varFISOUfig}). 
\par\noindent
Comparing the FIOU and FISOU cases of Figg.\ref{surfaceFIOUFISOU}-\ref{2dimFIOUFISOU}, the covariances 
exhibit a similar behavior for values of $\alpha$ greater than 0.8.
Furthermore, it appears evidente that the FIOU and FISOU  covariances 
are always lower than those of FIBM in Fig. \ref{2dimcolorcovFIBM}. The behaviors of  FIOU and FISOU covariance appear more similar each other, but really different from the behavior of the FIBM covariance. Indeed, in Fig.\ref{2dimcolorcovFIBM} is evident that the covariance of FIBM, as function of $\alpha$, and for increasing $\alpha$, attains rapidly values around 15, whereas the covariances of  FIOU and FISOU arrive to around 2, for $\alpha\approx 1.$

Moreover, referring to color plot of the covariances 
of FIOU and FISOU in Fig.\ref{2dimFIOUFISOU}, we note that
for every values of $\alpha$ the FISOU covariance, although slightly,  has values greater than those of the FIOU covariance.
In addition, the FISOU covariance show higher values (consequently, it can be seen a more diffuse correlation) around the diagonal respect to the case of FIOU covariance.
In particular, note that, for $\alpha = 0.2$ in Fig.\ref{2dimFIOUFISOU},
FISOU covariance $c(u,t)$ shows highest values for small values of $t$ and $u$ close to the diagonal, differently from the FIOU covariance. 
For $\alpha = 0.5$  and  $\alpha = 0.8$ the color plots of FISOU become more 
similar to those of FIOU, even if a more diffusion correlation remains evident for FISOU case. 
An explanation of these differences can be that the SOU process has a 
random starting point, and this has an impact on the ISOU, implying a greater variability 
(larger variance and correlation) respect to the case of IOU process.\par\noindent
The above comparing remarks can also be verified in  the three-dimensional plots in Fig.\ref{surfaceFIOUFISOU}, where it is also possible to
observe how the covariance functions increase for increasing values of $\alpha$.
\par\noindent
However, for large enough time $t,$  variance and covariance of the 
three processes are increasing functions of $\alpha,$ as one expects.

%

We remark the usefulness of the graphical comparisons for the considered processes and for different values of fractional order $\alpha$.
Again we confirm that the parameter $\alpha$ can be extremely useful in models, similar to the proposed neuronal  model, for tuning the length of memory and the accuracy of the time-scale of dynamics under observation.



The aim of giving further contribution to neuronal modeling  led us
to consider the long range memory process obtained as the fractional integral over time
of a GM process. In particular, the adoption of the fractional integral operator in place of the integer one, extending (form mathematical point of view) the latter, allows to describe dynamics that occur on different time-scales on which the integration of their past evolution plays the key rule to understand and predict the resultant complex behaviors.

In this paper, motivated by the above considerations, we  have studied as replacing the ordinary Riemann integral with the fractional RL
integral of order $\alpha \in (0,1)$ affects the behavior, when
varying $\alpha ,$ of an integrated GM process, i.e. the
processes obtained applying  the fractional RL integral over time to GM processes; in particular, we have studied
the fractional integral of Brownian motion,  non-stationary OU process and stationary OU process. These processes have a preeminent use in neuronal modeling and they are representative of the main classes of GM processes.
Beyond the closed-form obtained of covariance function of each of the above processes, we also provide numerical and graphical results, suitably compared.
The possible applications by means of a neuronal model are illustrated. The advantage of the evaluations of covariances for fractional integrated processes is in the possibility to implement a simulation procedure for sample paths and the consequent investigation of the first passage times through specific boundaries (\cite{Abu17},\cite{tai10}). This will be of extreme interest to explain and predict neuronal dynamics, because none theoretical result about these times for fractional integrated processes is available.
All these further investigations will be the object of our future work.

\end{document}